\input amstex


\def\Aut{\text{Aut}}
\def\b1{\text{\bf 1}}
\def\BA{{\Bbb A}}
\def\BC{{\Bbb C}}

\def\BP{{\Bbb P}}
\def\BZ{{\Bbb Z}}
\def\CA{{\Cal A}}
\def\CB{{\Cal B}}
\def\CD{{\Cal D}}

\def\CF{{\Cal F}}
\def\CL{{\Cal L}}
\def\CM{{\Cal M}}

\def\CO{{\Cal O}}
\def\CP{{\Cal P}}

\def\CT{{\Cal T}}
\def\CV{{\Cal V}}
\def\CW{{\Cal W}}

\def\dpar{\partial}
\def\End{\text{End}}

\def\id{\text{id}}
\def\fg{{\frak g}}

\def\gl{{\text{gl}}}
\def\Gr{{\text{Gr}}}
\def\hA{\widehat A}

\def\hCT{\widehat{{\Cal T}}}
\def\hfg{\widehat{\frak g}}
\def\hgl{\widehat{\text{gl}}}
\def\hOmega{\widehat{\Omega}}
\def\hsl{\widehat{\text{sl}}}
\def\hV{\widehat V}
\def\hX{\widehat X}
\def\Ker{\text{Ker}}

\def\hX{{\hat X}}
\def\Spec{\text{Spec}}
\def\Spf{\text{Spf}}
\def\ta{\tilde a}
\def\tb{\tilde b}
\def\tc{\tilde c}
\def\tCT{\tilde{{\Cal T}}}

\def\tG{\tilde G}
\def\tg{\tilde g}
\def\tI{\tilde I}
\def\tJ{\tilde J}
\def\tK{\tilde K}
\def\tL{\tilde L}
\def\tphi{\tilde\phi}
\def\tpsi{\tilde\psi}
\def\tQ{\tilde Q}
\def\Tr{\text{Tr}}

\def\tW{\tilde{W}}
\def\tX{\tilde{X}}


\def\btu{\bigtriangleup}
\def\hra{\hookrightarrow}
 
\def\lla{\longleftarrow}
\def\lra{\longrightarrow}
\def\twolra{\buildrel\longrightarrow\over\longrightarrow}

\parskip=6pt

\documentstyle{amsppt}
\document
\NoBlackBoxes


\centerline{\bf CHIRAL DE RHAM COMPLEX} 
\bigskip
\centerline{Fyodor Malikov, Vadim Schechtman and Arkady Vaintrob}
\bigskip

\bigskip

\bigskip
\centerline{\bf Introduction}
\bigskip

{\bf 0.1.} The aim of this note is to define certain sheaves of vertex 
algebras on smooth manifolds. In this note, "vertex algebra" will have 
the same meaning as in Kac's book [K]. Recall that these algebras are by 
definition $\BZ/(2)$-graded.  
"Smooth manifold" will mean a smooth scheme of finite type over $\BC$.    

For each 
smooth manifold $X$, we construct a sheaf 
$\Omega^{ch}_X$, called the {\bf chiral de Rham complex} of $X$. It is a sheaf 
of vertex algebras in the Zarisky topology, i.e. for each open $U\subset X$, 
$\Gamma(U;\Omega^{ch}_X)$ 
is a vertex algebra, and the restriction maps are morphisms of vertex 
algebras. It comes equipped with a $\BZ$-grading by {\it fermionic charge},
and the {\it chiral de Rham differential} $d_{DR}^{ch}$, which is    
an endomorphism of degree $1$ such that $(d_{DR}^{ch})^2=0$. One has 
a canonical embedding of the usual de Rham complex
$$
(\Omega_X, d_{DR})\hra (\Omega_X^{ch}, d_{DR}^{ch})
\eqno{(0.1)}
$$
The sheaf $\Omega^{ch}_X$ has also another $\BZ_{\geq 0}$-grading, by {\it 
conformal weight}, compatible 
with fermionic charge one. The differential $d^{ch}_{DR}$ respects conformal 
weight,  
and the subcomplex  
$\Omega_X$ coincides with the conformal weight zero component of 
$\Omega^{ch}_X$.   
The wedge multiplication on $\Omega_X$ may be restored from the 
operator product 
in $\Omega_X^{ch}$, see 1.3. 
The map (0.1) is a quasiisomorphism, cf. Theorem 4.4. Each component 
of $\Omega^{ch}_X$ of fixed conformal weight admits  
a canonical finite filtration whose graded factors are symmetric and 
exterior powers 
of the tangent bundle $\CT_X$ and of the bundle of $1$-forms $\Omega^1_X$.    

Similar sheaves exist in complex-analytic and $C^{\infty}$ settings. 

If $X$ is Calabi-Yau, then the sheaf $\Omega_X^{ch}$ has a 
structure of a {\it topological vertex algebra} (i.e. it admits $N=2$ 
supersymmetry, cf. 2.1), cf. 4.5. (For an arbitrary $X$, the obstruction 
to the existence of this structure is expressible in terms of the first  
Chern class of $\CT_X$, cf. Theorem 4.2.) 

One may hope that the vertex algebra $R\Gamma(X;\Omega^{ch}_X)$ defines 
the conformal field theory which is Witten's "$A$-model" associated 
with $X$.

The intuitive geometric picture behind our construction is as follows. 
Let $LX$ be the space of "formal loops" on $X$, i.e. of the 
maps of the punctured formal disk to $X$. Let $L^+X\subset LX$ be the 
subspace of loops regular at $0$. Note that we have a natural 
projection $L^+X\lra X$ (value at $0$). We have a functor 
$$
p:\ (Sheaves\ on\ LX)\lra (Sheaves\ on\ X), 
$$
namely, if $\CF$ is a sheaf on $LX$, then $\Gamma(U;p(\CF))=\Gamma(LU;F)$ 
for an open $U\subset X$.   
Now, the sheaf $\Omega^{ch}_X$ 
is the image under $p$ of the semiinfinite de Rham complex of the $\CD$-module 
of $\delta$-functions along $L^+X$. 

This sheaf is a particular case of a more general construction 
which associates with every $\CD$-module $\CM$ over $X$ its "chiral 
de Rham complex" $\Omega^{ch}_X(\CM)$ which is a sheaf 
of vertex modules over the vertex algebra $\Omega^{ch}_X$.  
Its construction is sketched in \S 6, cf. 6.11.       

{\bf 0.2.} One can also try to define a purely even sheaf 
$\CO_X^{ch}$ of vertex 
algebras, which could be called {\bf chiral structure sheaf}. Here the 
situation 
is more subtle than in the case of $\Omega^{ch}_X$, where  
"fermions cancel the anomaly". On can define this sheaf for curves, cf. 5.5.  
If $\dim(X)>1$, then there exists a non-trivial obstruction of 
cohomological nature to the construction of 
$\CO_X^{ch}$. This obstruction can 
be expressed in terms of a certain homotopy Lie algebra, cf. \S 5 A. 

However, one can define $\CO_X^{ch}$ for the flag manifolds $X=G/B$ ($G$ being 
a simple algebraic group and $B$ a Borel subgroup), cf. 
\S 5 B, C. This sheaf admits a structure of a $\hfg$-module at the 
critical level (here $\fg=Lie(G)$ and $\hfg$ is the corresponding 
affine Lie algebra).   
The space of global sections $\Gamma(X;\CO_X^{ch})$ is the irreducible vacuum 
$\hfg$-module for $\fg=sl(2)$; conjecturally, it is also true 
for any $\fg$, cf. 5.13. 
The sheaf $\CO_X^{ch}$ may be regarded as localization of Feigin-Frenkel 
bosonization.

More generally, if we start from an arbitrary $\CD$-module $\CM$ on $X=G/B$ 
correponding to some $\fg$-module $M$, then we can define its "chiralization" 
$\CM^{ch}$ which is an $\CO^{ch}$-module. It seems plausible that 
the space of global sections $\Gamma(X;\CM^{ch})$ coincides with the 
Weyl module over $\hfg$ corresponding to $M$ (on the critical level).           

{\bf 0.3.} A few words about the plan of the note. In \S 1 we recall the basic 
definitions, 
to fix terminology and notations. In \S 2, some "free field representation" 
results 
are described. No doubt, they are all well known (although for some of them 
we do not know the precise reference). They are a particular case of the  
construction of $\Omega_X^{ch}$, given in \S 3. 
In \S 4, we discuss the topological structure, and in \S 5 the chiral 
structure sheaf. In \S 6 we outline another construction of our 
vertex algebras, and some generalizations.   

{\bf 0.4.} The idea of this note arose from reading   
the papers [LVW] and [LZ]. We have learned the first example of 
"localization along the target space" from B.~Feigin, to whom goes 
our deep gratitude. We thank the referee for the useful remarks 
which helped to improve the exposition. We are thankful to 
V.~Gorbounov for catching many misprints.   

After the submission of this note to the alg-geom server, an interesting 
preprint [B] has appeared, where a possible application of the 
chiral de Rham complex to Mirror Symmetry was suggested.  

The work was done while the authors were visiting the Max-Planck-Institut 
f\"ur Mathematik in Bonn. We are grateful to MPI for the excellent working 
conditions, and especially to Yu.I.Manin for highly stimulating 
environment.  

\bigskip 

\bigskip
\centerline{\bf \S 1. Recollections on vertex algebras}
\bigskip

For more details on what follows, see [K]. For a 
coordinate-free exposition, see [BD1].  

{\bf 1.0.} For a vector space $V$, $V[[z,z^{-1}]]$ will 
denote the space of all formal 
sums $\sum_{i\in\BZ}\ a_iz^i,\ a_i\in V$. $V((z))$ will denote the 
subspace of Laurent power series, i.e. the sums as above, with 
$a_i=0$ for $i<<\infty$. We denote by $\dpar_z$ the operator of differentiation 
by $z$ acting on $V[[z,z^{-1}]]$. We set $\dpar_z^{(k)}:=(\dpar_z)^k/k!$, 
for an integer $k\geq 0$. For $a(z)\in V[[z,z^{-1}]]$, we will also use the 
notation $f(z)'$ for $\dpar_zf(z)$, and $\int a(z)$ for the coefficient at 
$z^{-1}$ (i.e. the residue).     

$V^*$ will denote the dual space. In the sequel, we will omit the prefix  
"super" in the words like "superalgebra", etc. 
For example,  "a Lie algebra" will mean  
"a Lie superalgebra".       

{\bf 1.1.} Let us define a {\bf vertex algebra} following [K], 4.1.   
Thus, a vertex algebra is the data $(V, Y, L_{-1}, \b1)$. Here $V$ is  
a $\BZ/(2)$-graded vector space $V=V^{ev}\oplus V^{odd}$. 
The $\BZ/(2)$-grading is called {\it parity};   
for an element $a\in V$, its parity will be denoted by $\ta$. The space of 
endomorphisms $\End(V)$ inherits the $\BZ/(2)$-grading from $V$.     

$\b1$ is an even element of $V$, called  
{\it vacuum}. $Y$ is an even linear mapping
$$
Y:\ V\lra\End(V)[[z,z^{-1}]]
\eqno{(1.1)}
$$
For $a\in V$, the power series $Y(a)$ will be denoted $a(z)$, and called 
the {\it field} corresponding to $a$. The coefficients of the series $a(z)$ 
are called {\it Fourier modes}.  

$L_{-1}$ is an even endomorphism of $V$.  
  
The following axioms must be satisfied.  

{\bf Vacuum Axiom.} $L_{-1}(\b1)=0; \b1(z)=\id$ (the constant power series). 
For all $a\in V$, $a(z)(\b1)\bigr|_{z=0}=a$ (in particular, the components 
of $a(z)$ at the negative powers of $z$ act as zero on the vacuum).   

{\bf Translation Invariance Axiom.} For all $a\in V$,  
$[L_{-1},a(z)]=\dpar_za(z)$

{\bf Locality Axiom.} For all $a, b\in V$, $(z-w)^N[a(z),b(w)]=0$ for 
$N>>0$. 

The meaning of this equality is explained in [K]. 

Given two fields $a(z), b(z)$, their {\it operator product expansion} is 
defined as in [K], (2.3.7a)
$$
a(z)b(w)=\sum_{j=1}^{N}\ \frac{c^j(w)}{(z-w)^j} + :a(z)b(w):
$$
where 
$$
:a(z)b(w):=a(z)_+b(w)+(-1)^{\ta\tb}b(w)a(z)_-
\eqno{(1.2)}
$$
Here 
$$
a(z)_-=\sum_{n\geq 0}\ a_{(n)}z^{-n-1};\ \ 
a(z)_+=\sum_{n<0}\ a_{(n)}z^{-n-1}
$$
for $a(z)=\sum\ a_{(n)}z^{-n-1}$, cf. [K], (2.3.5), (2.3.3).    

A {\it morphism} of vertex algebras $(V,\b1,L_{-1},Y)\lra 
(V',\b1',L'_{-1},Y')$ is  
an even linear map $f:\ V\lra V'$ taking $\b1$ to $\b1'$, such that 
$f\circ L_{-1}=L'_{-1}\circ f$, and  
for each $a\in V$, if $a_n$ is a Fourier mode of the 
field $a(z)$, we have $f(a)_n\circ f=f\circ a_n$. 

{\bf 1.2.} A {\bf conformal vertex algebra} (cf. [K], 4.10) is a vertex algebra 
$(V, Y, L_{-1}, \b1)$, together with an even element $L\in V$ such that 
if we write the field $L(z)$ as 
$$
L(z)=\sum_n\ L_nz^{-n-2},
$$
the endomorphisms $L_n$ satisfy the Virasoro commutation relations 
$$
[L_n,L_m]=(n-m)L_{n+m}+\frac{n^3-n}{12}\cdot c\cdot\delta_{n,-m}
\eqno{(1.3)}
$$
Here $c$ is a complex number, to be called the {\it Virasoro central charge} 
of $V$ (in [K] it is called the rank). 

The component $L_{-1}$ of the series (1.2) must coincide with the endomorphism 
$L_{-1}$ from the definition of a vertex algebra.  

The endomorphism $L_0$ must be diagonalizable, with integer eigenvalues. 
Set $V^{(n)}=\{ a\in V|\ L_0(a)=na\}$. We must have $V=\oplus_{n\in\BZ}\ 
V^{(n)}$. For $a\in V^{(n)}$, the number $n$ is called the {\bf conformal 
weight} of $a$ and denoted $|a|$. This $\BZ$-grading induces a 
$\BZ$-grading on the space $\End(V)$, also to be called the conformal 
weight. By definition, an endomorphism has conformal weight $n$ if it maps 
$V^{(i)}$ to $V^{(i+n)}$ for all $i$.  

The conformal weight grading on $V$ should be compatible 
with the parity grading, i.e. both gradings should come from a 
$\BZ\times\BZ/(2)$-bigrading. 

We require that if $a$ has conformal weight $n$, then the field $a(z)$ has 
the form
$$
a(z)=\sum_i\ a_{i}z^{-i-n},
\eqno{(1.4)}
$$
with $a_i$ having the conformal weight $-i$. 

We have the important {\bf Borcherds formula}. If $a, b$ are 
elements in a conformal vertex algebra, with $|a|=n$ then 
$$
a(z)b(w)=\sum_i\ \frac{a_i(b)(w)}{(z-w)^{i+n}}+:a(z)b(w):
\eqno{(1.5)}
$$
cf. [K], Theorem 4.6. 

The Virasoro commutation relations (1.3) are equivalent to the operator product 
$$
L(z)L(w)=\frac{c}{2(z-w)^4}+\frac{2L(w)}{(z-w)^2}+\frac{L(w)'}{z-w}
\eqno{(1.6)}
$$
(we omit from now on the regular part).         
  
{\bf 1.3.} Let $V$ be a conformal vertex algebra, such 
that $V^{(n)}=0$ for $n<0$.  
For $a\in V^{(n)},\ b\in V^{(m)}$, 
consider the operator product 
$$
a(z)b(w)=\sum_i\ \frac{a_i(b)(w)}{(z-w)^{i+n}}
\eqno{(1.7)} 
$$ 
We have $|a_i(b)|=m-i$. Therefore, $a_{m}(b)z^{-n-m}$ is the most 
singular term in (1.5).
 
If $n=m=0$ then $a(z)b(w)$ is non-singular at $z=w$. Define 
the multiplication on the space $V^{(0)}$ by the rule 
$$
a\cdot b=\biggl(a(z)b(w)\bigr|_{z=w}\b1\biggr)\biggr|_{z=0}
\eqno{(1.8)}
$$
This endows $V^{(0)}$ with a structure of an associative and commutative 
algebra with unity $\b1$.      

{\bf 1.4. Heisenberg algebra.} In this Subsection, we will define a conformal  
vertex algebra, to be called {\bf Heisenberg vertex algebra}. 

Fix an integer $N>0$. Let $H_N$ be the Lie algebra which as a vector space 
has the base $a^i_n,\ b^i_n,\ i=1,\ldots,N;\ n\in \BZ$, and 
$C$, all these elements being even, with the brackets
$$
[a^j_m,b^i_n]=\delta_{ij}\delta_{m,-n}C,
\eqno{(1.9)}
$$
all other brackets being zero. 

Our vertex algebra, to be denoted $V_N$, as a vector space is the 
{\bf vacuum representation} of $H_N$. As an $H_N$-module, it is 
generated by one vector \b1, subject to the relations  
$$
b^i_n \b1 =0
\ \ \text{if}\ n>0;\ 
a^i_m \b1 =0\ \ \text{if}\ m\geq 0;\ 
C \b1 = \b1 
\eqno{(1.10)}       
$$
The mapping 
$$
P(b^i_n,a^j_m)\mapsto P(b^i_n,a^j_m)\cdot \b1 
$$
identifies $V_N$ with the ring of commuting polynomials on the variables 
$b^i_n,\ a^j_m,\ n\leq 0,\ m<0,\ i,j=1,\ldots,N$. We will use 
this identification below.   

Let us define the structure of a vertex algebra on $V_N$. The $\BZ/(2)$ grading 
is trivial: everything is even. The vacuum vector is \b1. 

The fields corresponding to the elements of $V_N$ are defined by induction 
on the degree of a monomial. First we define the fields for the degree 
one monomials by setting
$$
b^i_0(z)=\sum_{n\in\BZ}\ b^i_nz^{-n};\ 
a^j_{-1}(z)=\sum_{n\in\BZ}\ a^j_nz^{-n-1}
\eqno{(1.11)}
$$
Here $b^i_n,\ a^j_n$ in the right hand side are regarded as 
operating on $V_N$ by multiplication. 

We set 
$$
b^i_{-n}(z)=\dpar_z^{(n)}b^i_0(z),\ \ 
a^j_{-n-1}(z)=\dpar_z^{(n)}a^i_{-1}(z)
\eqno{(1.12)}
$$
for $n>0$. 

The fields corresponding to the monomials of degree $>1$ are defined by 
induction, using the {\it normal ordering}. Let us call the 
operators $b^i_n,\ n>0$, and $a^j_n,\ n\geq 0$, acting on 
$V_N$, {\it annihilation 
operators}. 

For $x=b^i_n$ or $a^j_n,\ n\in\BZ$, and 
$b\in \End(V_N)$, the {\bf normal ordered product} 
$:xb:$ is given by 
$$
:xb:=\left\{\aligned bx&\text{\ \ \ if $x$ is an annihilation operator}\\
xb&\text{\ \ \ otherwise}\endaligned\right. 
\eqno{(1.13)}
$$
Define by induction 
$$
:x_1\cdot\ldots\cdot x_k:=:x_1\cdot (:x_2\cdot\ldots\cdot x_k:):, 
\eqno{(1.14)}
$$
for $x_p=b^i_n$ or $a^j_n$, $p=1,\ldots, k$. 

Given two series $x(z)=\sum_{n\in\BZ}\ x_nz^{-n+p}$ and 
$y(z)=\sum_{n\in\BZ}\ y_nz^{-n+q}$, with $x_n$ as above, we set
$$
:x(z)y(w):=\sum_{n,m\in\BZ}\ :x_ny_m:\ z^{-n+p}w^{-m+q}
\eqno{(1.15)}
$$
For any finite sequence $x_1,\ldots,x_p$, where each $x_j$ is equal to 
one of $a^i_n$ or $b^i_n$, we define the series 
$:x_1(z)x_2(z)\cdot\ldots\cdot x_p(z):\in\End(V_N)[[z,z^{-1}]]$ 
by induction, as in (1.14). This expression does not depend on the 
order of $x_i$'s.  

Given a monomial $x_1\cdot\ldots\cdot x_p\b1\in V_N$, with $x_i$ as above, 
we define the corresponding field by 
$$
x_1\cdot\ldots\cdot x_p(z)=:x_1(z)\cdot\ldots\cdot x_p(z):
\eqno{(1.16)}
$$
Since every element of $V_N$ is a finite linear combination of 
monomials as above,   
this completes the definition of the mapping (1.1).
 
We will use the shorthand notations
$$
b^i(z)=b^i_0(z);\ a^j(z)=a^j_{-1}(z)
\eqno{(1.17)}
$$
The operator products of these basic fields are 
$$
a^j(z)b^i(w)=\frac{\delta_{ij}}{z-w}+\text{(regular)}
\eqno{(1.18a)}
$$
$$
b^i(z)b^j(w)=\text{(regular)};\ a^i(z)a^j(w)=\text{(regular)}
\eqno{(1.18b)}
$$
where "(regular)" means the part regular at $z=w$. 
These operator products are equivalent to the commutation relations (1.9).  

Other operator products are computed by differentiation of (1.18), 
and using Wick theorem, cf. [K], 3.3. 

One can say that the vertex algebra $V_N$ is generated by the even  
fields $b^i(z),\ a^j(z)$, of conformal weights 
$0$ and $1$ respectively, subject to the relations 
(1.18). 
  
The Virasoro field is given by 
$$
L(z)=\sum_{i=1}^N\ :b^i(z)'\cdot a^i(z):
\eqno{(1.19)}
$$
The central charge is equal to $2N$. Let us check this. Assume for simplicity 
that $N=1$, and let us omit the index $1$ at the fields $a, b$. Thus, we have 
$$
L(z)=:b(z)'a(z):
$$
Let us compute the operator product $L(z)L(w)$ using the Wick theorem. We have 
$$
:b(z)'a(z)::b(w)'a(w):=[b(z)'a(w)][a(z)b(w)']+
$$
$$
+[b(z)'a(w)]:a(z)b(w)':+[a(z)b(w)']:b(z)'a(w):=
$$
(we have $b(z)'a(w)=a(z)b(w)'=1/(z-w)^2$)
$$
=\frac{1}{(z-w)^4}+\frac{2:b(w)'a(w):}{(z-w)^2}+\frac{:b(w)''a(w):+
:b(w)'a(w)':}{z-w}
$$
Hence,
$$
L(z)L(w)=\frac{1}{(z-w)^4}+\frac{2L(w)}{(z-w)^2}+\frac{L(w)'}{z-w}
$$
which is (1.6) with $c=2$.   

{\bf 1.5. Clifford algebra.} Let $Cl_N$ be the Lie algebra which as a vector 
space has the base $\phi^i_n,\ \psi^i_n,\ i=1,\ldots,N,\ n\in \BZ$, and $C$, 
all these elements being {\it odd}, with the brackets
$$
[\phi^j_m,\psi^i_n]=\delta_{ij}\delta_{m,-n}\cdot C
\eqno{(1.20)}
$$
{\bf Clifford vertex algebra} $\Lambda_N$ is defined as in the previous 
subsection, starting with the odd fields
$$
\psi^i(z)=\sum_{n\in\BZ}\ \psi^i_nz^{-n-1};\ 
\phi^j(z)=\sum_{n\in\BZ}\ \phi^j_nz^{-n}
\eqno{(1.21)}
$$ 
and repeating the 
definitions of {\it loc. cit.}, with  
$a$ (resp. $b$) replaced by $\psi$ (resp. $\phi$).  
One must put the obvious signs 
in the definition of the normal ordering. Thus, $\Lambda_N$ is generated 
by the odd fields $\phi^i(z), \psi^i(z)$, subject to the relations
$$
\phi^i(z)\psi^j(w)=\frac{\delta_{ij}}{z-w}+\text{regular};
\eqno{(1.22a)}
$$
$$
\phi^i(z)\phi^j(w)=\text{regular};\ \psi^i(z)\psi^j(w)=\text{regular}
\eqno{(1.22b)}
$$ 
The Virasoro field is given by 
$$
L(z)=\sum_{i=1}^N\ :\phi^i(z)'\cdot\psi^i(z):
\eqno{(1.23)}
$$
The central charge is equal to $-2N$. 

{\bf 1.6.} If $\CA$ and $\CB$ are vertex algebras then their tensor product 
$\CA\otimes\CB$ admits a canonical structure of a vertex algebra, cf. [K], 
4.3. The Virasoro element is given by 
$$
L_{\CA\otimes\CB}=L_{\CA}\otimes 1+1\otimes L_{\CB}
\eqno{(1.24)}
$$ 
We will use in the sequel the tensor product of the Heisenberg and Clifford 
vertex algebras $\Omega_N=V_N\otimes\Lambda_N$. Its Virasoro central charge  
is equal to $0$.

{\bf 1.7.} Let $\CA$ be a vertex algebra. A linear map $f:\ \CA\lra\CA$ is 
called {\it derivation} of $\CA$ if for any $a\in\CA$, 
$$
f(a)(z)=[f,a(z)]
\eqno{(1.25)}
$$
Note that an invertible map $g:\ \CA\lra\CA$ is an automorphism of $\CA$ iff 
$$
g(a)(z)=ga(z)g^{-1}
\eqno{(1.26)}
$$
and (1.25) is the infinitesimal version of (1.26). 

It follows from Borcherds formula that for every $a\in\CA$, the Fourier 
mode $\int a(z)$ is a derivation, cf. Lemma 1.3 from [LZ]. 
Consequently, if the endomorphism 
$\exp(\int a(z))$ is well defined, it is an automorphism of $\CA$.   

\bigskip 

\bigskip
\centerline{\bf \S 2. De Rham chiral algebra of an affine space}
\bigskip

{\bf 2.1.} A {\bf topological vertex algebra of rank d} is 
a conformal vertex algebra $\CA$ of Virasoro central charge $0$,  
equipped with an even element $J$ of conformal weight $1$, and  
two odd elements, $Q$, of conformal weight $1$, 
and $G$, of conformal weight $2$.  

The following operator products must hold  
$$
L(z)L(w)=\frac{2L(w)}{(z-w)^2}+\frac{L(w)'}{z-w}
\eqno{(2.1a)}
$$
$$  
J(z)J(w)=\frac{d}{(z-w)^2};\ 
L(z)J(w)=-\frac{d}{(z-w)^3}+\frac{J(w)}{(z-w)^2}+\frac{J(w)'}{z-w}
\eqno{(2.1b)}
$$
$$
G(z)G(w)=0;\ L(z)G(w)=\frac{2G(w)}{(z-w)^2}+\frac{G(w)'}{z-w};\ 
J(z)G(w)=-\frac{G(w)}{z-w}
\eqno{(2.1c)}
$$
$$
Q(z)Q(w)=0;\ L(z)Q(w)=\frac{Q(w)}{(z-w)^2}+\frac{Q(w)'}{z-w};\ 
J(z)Q(w)=\frac{Q(w)}{z-w}
\eqno{(2.1d)}
$$
$$
Q(z)G(w)=\frac{d}{(z-w)^3}+\frac{J(w)}{(z-w)^2}+\frac{L(w)}{z-w}
\eqno{(2.1e)}
$$
Note the following consequence of (2.1e), 
$$
[Q_0,G(z)]=L(z)
\eqno{(2.2)}
$$ 

{\bf 2.2.} In this subsection, 
we will intorduce a structure of a topological vertex algebra of rank $N$ 
on the vertex algebra $\Omega_N$ from 1.6. This topological vertex algebra 
will be called {\bf de Rham chiral algebra of the 
affine space $\BA^N$}. 

Recall that the Virasoro element is given by
$$
L=\sum_{i=1}^N\ \bigl(b_{-1}^ia^i_{-1}+\phi^i_{-1}\psi^i_{-1}\bigr)
\eqno{(2.3a)}
$$ 
Define the elements $J, Q, G$ by 
$$
J=\sum_{i=1}^N\ \phi_0^i\psi_{-1}^i;\ 
Q=\sum_{i=1}^N\ a_{-1}^i\phi_0^i;\ 
G=\sum_{i=1}^N\ \psi_{-1}^ib_{-1}^i
\eqno{(2.3b)}
$$
The corresponding fields are
$$
L(z)=\sum\ \bigl(:b^i(z)'a^i(z):+:\phi^i(z)'\psi^i(z):\bigr)
\eqno{(2.4a)}
$$
and  
$$
J(z)=\sum\ :\phi^i(z)\psi^i(z):,\ 
Q(z)=\sum\ :a^i(z)\phi^i(z):,\ 
G(z)=\sum\ :\psi^i(z)b^i(z)':
\eqno{(2.4b)}
$$
The relations (2.1) are readily checked using the Wick theorem. 

{\bf 2.3.} Let us define the {\bf fermionic charge} operator acting on 
$\Omega_N$, by 
$$
F:=J_0=\sum_i\sum_n\ :\phi^i_n\psi^i_{-n}:
\eqno{(2.5)}
$$
We have
$$
F\b1=0
\eqno{(2.6)}
$$
and
$$
[F,\phi^i_n]=\phi^i_n;\ [F,\psi^i_n]=-\psi^i_n;\ [F,a^i_n]=[F,b^i_n]=0
\eqno{(2.7)}
$$
We set
$$
\Omega^p_N=\{\omega\in\Omega_N|\ F\omega=p\omega\}
\eqno{(2.8)}
$$
Obviously,
$$
\Omega_N=\oplus_{p\in\BZ}\ \Omega^p_N
\eqno{(2.9)}
$$
We define an endomorphism $d$ of the space $\Omega_N$ by  
$$
d:=-Q_0=-\sum_{i, n}\ :a^i_n\phi^i_{-n}:
\eqno{(2.10)}
$$
(we could omit the normal ordering in the last formula since the letters 
$a$ and $\phi$ commute anyway). 
We have $d^2=0$. Indeed, by Wick theorem $Q(z)Q(w)=\text{regular}$, hence 
all Fourier modes of $Q(z)$ (anti)commute. The map $d$ is 
called {\bf chiral de Rham differential}.
  
The map $d$ increases the fermionic charge by $1$, by (2.7). 
Thus, the space $\Omega_N$ equipped with the fermionic charge grading 
and the differential $d$, becomes a complex (infinite in both directions), 
called {\bf chiral de Rham complex} of $\BA^N$.   

Consider the usual algebraic de Rham complex $\Omega(\BA^N)=\oplus_{p=0}^N\ 
\Omega^p(\BA^N)$ of 
the affine space $\BA^N$. We identify the coordinate functions 
with the letters $b_0^{1},\ldots,b_0^{N}$, and their differentials 
with the fermionic variables $\phi_0^1,\ldots,\phi_0^N$. 

Thus, we identify the commutative dg algebras 
$$
\Omega(\BA^N)=\BC[b_0^{1},\ldots,b_0^{N}]\otimes
\Lambda(\phi_0^1,\ldots,\phi_0^N),
\eqno{(2.12)}
$$
the second factor being the exterior algebra. The grading is 
defined by assigning to the letters $b_0^{i}$ (resp. $\phi_0^j$) the degree 
$0$ (resp. $1$). 

The usual de Rham differential is given by 
$$
d_{DR}=\sum_i\ a_0^i\phi_0^i,
\eqno{(2.13)}
$$
as follows from the relations (1.9).

{\bf 2.4. Theorem.} {\it The obvious embedding of complexes  
$$
i:\ (\Omega(\BA^N),d_{DR})\lra (\Omega_N,d)
\eqno{(2.14)}
$$
is compatible with the differentials, and is a quasiisomorphism.}

We identify the space $\Omega_N$ with the space of polynomials in the 
letters $b^i_n, \phi^i_n\ (n\leq 0)$ and $a^i_n, \psi^i_n\ (n<0)$. One sees 
that on the subspace $\BC[b^i_0,\phi^i_0]$, all the summands $a^i_n\phi^i_{-n}$ 
with $n\neq 0$ act trivially. It follows that the map $i$ is compatible 
with the differentials. 

To prove that it is a quasiisomorphism, let us split $d$ in two 
commuting summands $d=d_++d_-$ where
$$
d_+=\sum_i\sum_{n\geq 0}\ a^i_n\phi^i_{-n},\ 
d_-=\sum_i\sum_{n<0}\ a^i_n\phi^i_{-n}
\eqno{(2.15)}
$$
We think of the space $\Omega_N$ as of the tensor product 
$\BC[a^i_n,\psi^i_n]\otimes\BC[b^j_m,\phi^j_m]$. The differential $d_-$ 
acts trivially on the second factor, and on the first one it is the Koszul 
differential. So, the cohomology of $d_-$ is 
$\BC[b^j_{m},\phi^j_{m}]_{m\leq 0}$.
 
Now, we have to compute the cohomology of this space with respect to $d_+$. 
For this purpose, split $d_+$ once again as $d_+=d_{DR}+d_+^{>0}$, 
and our space as $\BC[b^j_0,\phi^j_0]\otimes\BC[b^j_m,\phi^j_m]_{m<0}$. 
The differential $d_+^{>0}$ acts trivially on the first factor, 
and on the second one, it is the de Rham differential. Hence (by Poincar\'e 
lemma), taking the cohomology of $d_+^{>0}$ 
kills all non-zero modes, and we are left precisely with the usual 
de Rham complex.

Alternatively, it follows from (2.2) that 
$$
[G_0,d]=L_0
\eqno{(2.16)}
$$
The operator $G_0$ commutes with $L_0$, and it follows from (2.16) 
that it gives a homotopy to $0$ for the operator $d$ on all 
the subcomplexes of non-zero conformal weight. Therefore, all cohomology 
lives in the conformal weight zero subspace. $\btu$

{\bf 2.5.} The vertex algebra $\Omega_N$ satisfies the assumptions 
of 1.3. The subspace $\Omega(\BA^N)$ coincides with the conformal 
weight zero component of it. If we apply the definition of 1.3, we get 
the structure of a commutative algebra on $\Omega(\BA^N)$ which is given  
by the usual wedge product of differential forms.  

\bigskip

\bigskip
\centerline{\bf \S 3. Localization}
\bigskip

{\bf 3.1.} Consider the Heisenberg vertex algebra $V_N$ defined in 1.4. 
As in {\it loc. cit.}, we will identify the space $V_N$ with the 
space of polymonials $\BC[b^i_{-n},a^j_{-m}]_{n\geq 0,\ m>0}$. To simplify 
the notations below, let us denote the zero mode variables $b^i_0$ 
by $b^i$. 

Let $A_N$ denote the algebra of polynomials $\BC[b^1,\ldots,b^N]$. 
The space $A_N$ is identified with the subspace 
$V_N^{(0)}\subset V_N$ of conformal weight zero.    
The space $V_N$ has an obvious structure of an $A_N$-module. 
Let $\hA_N$ denote the algebra of formal power series $\BC[[b^1,\ldots,b_N]]$. 
Set
$$
\hV_N=\hA_N\otimes_{A_N}V_N
\eqno{(3.1)}
$$
We are going to introduce a structure of a conformal vertex algebra 
on the space $\hV_N$. Let us define the map (1.1). 

Let $f(b^1,\ldots,b^N)$ be a power series from $\hA_N$. We claim that 
the expression $f(b^1(z),\ldots,b^N(z))$ makes sense as an element 
of $\End(\hV_N)[[z,z^{-1}]]$. (We are grateful to Boris Feigin 
who has shown us a particular case of the following construction.) 
Let us express the power series $b^i(z)$ as 
$$
b^i(z)=b^i+\Delta b^i(z)
\eqno{(3.2)}
$$
Thus, 
$$
\Delta b^i(z)=\sum_{n>0}\ (b^i_nz^{-n}+b^i_{-n}z^n)
\eqno{(3.3)}
$$
Let us define $f(b^1(z),\ldots,b^N(z))$ by the Taylor formula
$$
f(b^1(z),\ldots,b^N(z))=\sum\ \Delta b^1(z)^{i_1}\cdot\ldots\cdot \Delta 
b^N(z)^{i_N}
\dpar^{(i_1,\ldots,i_N)}f(b^1,\ldots,b^N)
\eqno{(3.4)} 
$$
where
$$
\dpar^{(i_1,\ldots,i_N)}=\frac{\dpar^{i_1}_{b^1}}{i_1!}\cdot\ldots\cdot
\frac{\dpar^{i_N}_{b^N}}{i_N!}
\eqno{(3.5)}
$$
We will show that the series (3.4) gives a well-defined element of 
$\End(\hV_N)[[z,z^{-1}]]$. 
Let us write
$$
\Delta b^1(z)^{i_1}\cdot\ldots\cdot \Delta b^N(z)^{i_N}=\sum_k\ 
c_k^{i_1,\ldots,i_N}z^{-k}
$$
The coefficient $c_k^{i_1,\ldots,i_N}$ is an infinite sum of the 
monomials 
$$
b_{k_1}^{j_1}\cdot\ldots\cdot b_{k_I}^{j_I},
\eqno{(3.6)} 
$$
with $I=i_1+\ldots+i_N$, $k_1+\ldots+k_I=k$. Pick an element $v\in \hV^N$. 
There exists $M$ such that $b_{l_1}^{j_1}\cdot\ldots\cdot b_{l_N}^{j_N}v=0$ 
if $\sum\ l_i>M$. 

We have 
$$
k=\sum\ k_i=\ ^+\sum |k_i|-\ ^-\sum\ |k_i|
$$
where $^+\sum$ (resp. $^-\sum$) denotes the sum of all positive (resp. 
negative) summands. If $b_{k_1}^{j_1}\cdot\ldots\cdot b_{k_I}^{j_I}v\neq 0$,  
then 
$$
^+\sum\ |k_i|\leq M
\eqno{(3.7a)}
$$
On the other hand,
$$
^-\sum\ |k_i|=\ ^+\sum\ |k_i|-k\leq M-k
\eqno{(3.7b)}
$$
There exists only a finite number of tuples $(k_1,\ldots,k_I)$ satisfying 
(3.7a) and (3.7b). Therefore, $c_k^{i_1,\ldots,i_N}$ are well-defined 
endomorphisms of $\hV_N$. 

We have
$$
f(b^1(z),\ldots,b^N(z))=\sum_k\ \biggl(\sum_{i_1,\ldots,i_N}
\ c_k^{i_1,\ldots,i_N}\dpar^{(i_1,\ldots,i_N)}f(b^1,\ldots,b^N)\biggr)
\ z^{-k}
\eqno{(3.8)}
$$
All numbers $k_i$ are non-zero. Let $I_+$ (resp., $I_-$) be the number of 
positive (resp. negative $k_i$'s). We have $I_-\leq\ ^-\sum k_i\leq M-k$, 
hence 
$M\geq\ ^+\sum\ k_i\geq I_+=I-I_-\geq I-M+k$. Therefore, 
$$
I\leq 2M-k
\eqno{(3.9)}
$$
Therefore, when we apply the series (3.8) to the element $v$, only a finite 
number of terms in the sum over $(i_1,\ldots,i_N)$ survives. Therefore, 
the series (3.8) is a well-defined element of $\End(\hV_N)[[z,z^{-1}]]$. 

Every element of $\hV_N$ is a finite sum of products 
$g(a)f(b)$ where $g(a)$ is a polynomial in the letters $a$ and $f(b)$ 
is a power series as above. We have already defined $f(b)(z)$. 
The definition of $g(a)(z)$ is the same as in the case of $V_N$. 
We define $\bigl(g(a)f(b)\bigr)(z)$ by 
$$
\bigl(g(a)f(b)\bigr)(z)=:g(a)(z)f(b)(z):
\eqno{(3.10)}
$$
where the normal ordering is defined in (1.2). This completes the definition 
of the mapping 
$$
Y:\ \hV_N\lra\End(\hV_N)[[z,z^{-1}]]
\eqno{(3.11)}
$$
The following version of the definition of the map (3.11) is helpful 
in practice. Every element $c\in\hV_N$ is a limit of the elements of 
$c_i\in V_N$ (in the obvious topology). We can regard the fields 
$c_i(z)$ as the elements of $\End(\hV_N)[[z,z^{-1}]]$. The field 
$c(z)$ is the limit of the fields $c_i(z)$.  

We define the vacuum and Virasoro element $\b1, L\in \hV_N$ as the image 
of the corresponding element of $V_N$ under the natural map $V_N\lra\hV_N$. 

{\bf 3.2. Theorem.} {\it The construction of} 3.1 {\it defines a structure 
of a conformal vertex algebra on the space $\hV_N$.}

This follows from [K], Theorem 4.5.  

{\bf 3.3.} Let $A_N^{an}\subset \hA_N$ denote the subalgebra of power 
series, convergent in a neighbourhood of the origin. Set 
$$
V_N^{an}=A_N^{an}\otimes_{A_N}V_N\subset \hV_N
\eqno{(3.12)}
$$
It is clear from the inspection of the Taylor formula (3.4) 
that for $f(a,b)\in V_N^{an}$, the Fourier modes of $f(a,b)(z)$ which 
belong to $\End(\hV_N)$, respect the subspace $V_N^{an}$. Therefore, 
the conformal vertex algebra structure on $\hV_N$ defined in 3.1 
induces the structure of a conformal vertex algebra on $V_N^{an}$. 

In this argument, $A_N^{an}$ can be replaced by any algebra of functions 
containing $A_N$ and closed under derivations.  
More precisely, one has the following general statement. Let $A'$  be an 
arbitrary commutative $A_N$-algebra, given together with an action 
of the Lie algebra $\CT=Der(A_N)$ by derivations, extending the natural action 
of $\CT$ on $A_N$. Then the space $V_{A'}:=A'\otimes_{A_N}V_N$ admits a natural 
structure of a vertex algebra. For the details, see 6.9.   

For example, let $A_N^{sm}$ denote the algebra of germs of smooth ($C^{\infty}$) 
functions.
Then we get a vertex algebra  
$$
V_N^{sm}=A_N^{sm}\otimes_{A_N}V_N
\eqno{(3.13)}
$$

Another natural example is that of localization of $A$. It is treated in the 
next 
subsection.  

{\bf 3.4. Zariski localization.} Let $f\in A_N$ be a nonzero polynomial. Let 
$A_{N;f}$ 
denote the localization $A_N[f^{-1}]$. Set 
$$
V_{N;f}=A_{N;f}\otimes_{A_N}V_N
\eqno{(3.14)}
$$
Consider the Taylor formula (3.4) applied to the function $f^{-1}$. 
We have evidently 
$$
\biggl[\dpar^{(i_1,\ldots,i_N)}f^{-1}\biggr](b^1,\ldots,b^N)\in A_{N;f}
$$

In more concrete terms, let $f(z)=\sum f_nz^{-n}$ is the field correponding to 
$f$, then 
we want to define the field corresponding to $f^{-1}$ as 
$$ 
f(z)^{-1}=(f_0+f_{-1}z+f_1z^{-1}+\ldots)^{-1}=
f_0^{-1}\bigl(1+f_0^{-1}(f_{-1}z+f_1z+\ldots)\bigr)^{-1}=
$$
(we use the geometric series)
$$
=f_0^{-1}(1+f_0^{-2}(2f_{-1}f_1+2f_{-2}f_2+\ldots)+\ldots)
$$
(we started to write down the coefficient at $z^0$). Now, in the right hand 
side, 
the coefficient at each power of $z$ is an infinite sum, but as an operator 
acting on 
$A_{N;f}$ it is well defined since only finite number of terms act nontrivially. 
We need only to invert $f_0=f$. 

Therefore, the construction 3.1 provides a conformal 
vertex algebra structure on the space $V_{N;f}$. 

Let $X$ denote the affine space $\Spec(A_N)$. Let $\CO_X^{ch}$ be the 
$\CO_X$-quasicoherent sheaf corresponding to the $A_N$-module $V_N$. 
We have just defined the structure of a conformal vertex algebra on the 
spaces $V_{N;f}=\Gamma(U_f;\CO_X^{ch})$ where $U_f=\Spec(A_{N;f})$. If 
$U_f\subset U_g$ then the restriction map $V_{N;g}\lra V_{N;f}$ is a 
morphism of conformal vertex algebras.  
 
If $U\subset X$ is an arbitrary open, we have $U=\bigcup U_f$, and 
$$
\Gamma(U;\CO_X^{ch})=\Ker\bigl(\prod\ V_{N;f}\twolra\prod\ V_{N;fg}\bigr)
$$
Using this formula, we get a structure of a conformal vertex algebra 
on the space $\Gamma(U;\CO_X^{ch})$. Therefore, $\CO_X^{ch}$ gets a structure 
of a sheaf of conformal vertex algebras. 

{\bf 3.5.} We can add fermions to the picture. Consider the spaces 
$\hOmega_N=\hA_N\otimes_{A_N}\Omega_N$, 
$\Omega_N^{an}=A_N^{an}\otimes_{A_N}\Omega_N\subset\hOmega_N$, 
$\Omega_N^{sm}=A_N^{sm}\otimes_{A_N}\Omega_N$.  
The construction 3.1 provides a structure of a topological 
vertex algebras on these spaces. 

Let $X$ be as in 3.4; let $\Omega_X^{ch}$ denote the $\CO_X$-quasicoherent 
sheaf associated with the $A_N$-module $\Omega_N$. The construction 3.4 
provides a structure of a sheaf of topological vertex algebras of rank $N$ 
on $\Omega_X^{ch}$.

{\bf 3.6.} Now we want to study coordinate changes in our vertex algebras. 
Let $X$ be the formal scheme $\Spf(\BC[[b^1,\ldots,b^N]])$. Consider the 
formal $N|N$-dimensional superscheme $\tX=\Pi TX$ (here $TX$ is the total 
space of the tangent bundle, $\Pi$ is the parity change functor). Thus, 
$\tX$ has the same underlying space as $X$, and the structure sheaf of 
$\tX$ coincides with the de Rham algebra of differential forms $X$. On $\tX$, we 
have $N$  
even coordinates $b^1,\ldots,b^N$ and odd ones $\phi^1=db^1,\ldots,
\phi^N=db^N$. 

To this superscheme $\tX$, with the above coordinates, we have assigned a 
(super)vertex algebra $\hOmega_N$, generated by the fields 
$b^i(z), a^i(z)$ (even ones) and $\phi^i(z), \psi^i(z)$ (odd ones). The fields 
$a^i(z)$ (resp. $\psi^i(z)$) correspond to the vector fields 
$\dpar_{b^i}$ (resp. $\dpar_{\phi^i}$) on $\tX$. 

These fields satisfy the relations (cf. (1.18), (1.22))
$$
a^i(z)b^j(w)=\frac{\delta_{ij}}{z-w}
\eqno{(3.15a)}
$$
$$
b^i(z)b^j(w)=\text{(regular)};\ a^i(z)a^j(w)=\text{(regular)}
\eqno{(3.15b)}
$$
$$
\phi^i(z)\psi^j(w)=\frac{\delta_{ij}}{z-w}
\eqno{(3.15c)}
$$
$$
\phi^i(z)\phi^j(w)=\text{(regular)};\ \psi^i(z)\psi^j(w)=\text{(regular)}
\eqno{(3.15d)}
$$
$$
b^i(z)\phi^j(w)=\text{(regular)};\  b^i(z)\psi^j(w)=\text{(regular)}
\eqno{(3.15e)}
$$
$$
a^i(z)\phi^j(w)=\text{(regular)};\ a^i(z)\psi^j(w)=\text{(regular)}
\eqno{(3.15f)}
$$
Consider an invertible coordinate transformation on $X$, 
$$
\tb^i=g^i(b^1,\ldots,b^N);\ b^i=f^i(\tb^1,\ldots,\tb^N)
\eqno{(3.16a)}
$$
where $g^i\in\BC[[b^j]];\ f^i\in\BC[[\tb^j]]$. It induces the transformation 
of the odd coordinates $\phi^i=db^i$, 
$$
\tphi^i=\frac{\dpar g^i}{\dpar b^j}\phi^j;\ 
\phi^i=\frac{\dpar f^i}{\dpar\tb^j}\tphi^j
\eqno{(3.16b)}
$$
(the summation over the repeating indices is tacitly assumed). 
The vector fields transform as follows, 
$$
\dpar_{\tb^i}=\frac{\dpar f^j}{\dpar\tb^i}(g(b))\dpar_{b^j}+
\frac{\dpar^2f^k}{\dpar\tb^i\dpar\tb^l}(g(b))\cdot\frac{\dpar g^l}{\dpar b^r}
\cdot\phi^r\dpar_{\phi^k}
\eqno{(3.16c)}
$$
and
$$ 
\dpar_{\tphi^i}=\frac{\dpar f^j}{\dpar\tb^i}(g(b))\dpar_{\phi^j}
\eqno{(3.16d)}
$$
We want to lift the transformation (3.16a) to the algebra $\Omega_N$. 
Define the tilded fields by 
$$
\tb^i(z)=g^i(b)(z)
\eqno{(3.17a)}
$$
$$
\tphi^i(z)=\biggl(\frac{\dpar g^i}{\dpar b^j}\phi^j\biggr)(z)
\eqno{(3.17b)}
$$
$$
\ta^i(z)=\biggl(a^j\frac{\dpar f^j}{\dpar\tb^i}(g(b))\biggr)(z)+
\biggl(\frac{\dpar^2 f^k}{\dpar\tb^i\dpar\tb^l}(g(b))
\frac{\dpar g^l}{\dpar b^r}\phi^r\psi^k\biggr)(z)
\eqno{(3.17c)}
$$
$$
\tpsi^i(z)=\biggl(\frac{\dpar f^j}{\dpar\tb^i}(g(b))\psi^j\biggr)(z)
\eqno{(3.17d)}
$$
 
{\bf 3.7. Theorem.} {\it The fields $\tb^i(z), \ta^i(z), \tphi^i(z)$ 
and $\tpsi^i(z)$ satisfy the relations} (3.15). 

We will use the relations
$$
h(b)(z)a^i(w)=-\frac{\dpar h/\dpar b^i(w)}{z-w};\ 
a^i(z)h(b)(w)=\frac{\dpar h/\dpar b^i(w)}{z-w}  
\eqno{(3.18)}
$$
for each $h\in\hA_N$, which follow from (3.15a) and the Wick theorem. 
Let us check (3.15a) for the tilded fields. We have 
$$
\tb^i(z)\ta^j(w)=g^i(b)(z)a^k\frac{\dpar f^k}{\dpar\tb^j}(g(b))(w)-
$$
$$
-g^i(b)(z)\frac{\dpar^2f^k}{\dpar\tb^j\dpar\tb^l}
(g(b))\frac{\dpar g^l}{\dpar b^r}\psi^k\phi^r(w)
$$
By the Wick theorem, 
the first summand is equal to
$$
-\frac{\dpar g^i}{\dpar b^k}\frac{\dpar f^k}{\dpar\tb^j}(g(b))(w)
\cdot\frac{1}{z-w}=
\frac{\delta_{ij}}{z-w}, 
$$
by (3.18). The second summand is zero. Let us check (3.15b). The first 
identity is clear. We have
$$
\ta^i(z)\ta^j(w)=a^k\frac{\dpar f^k}{\dpar\tb^i}(g(b))(z)
a^n\frac{\dpar f^n}{\dpar\tb^j}(g(b))(w)
+a^k\frac{\dpar f^k}{\dpar\tb^i}(g(b))(z)
\frac{\dpar^2f^n}
{\dpar\tb^j\dpar\tb^l}(g(b))\frac{\dpar g^l}{\dpar b^r}\psi^n\phi^r(w)-
$$
$$
-\frac{\dpar^2f^k}{\dpar\tb^i\dpar\tb^l}(g(b))\frac{\dpar g^l}{\dpar b^r}
\psi^k\phi^r(z)
a^n\frac{\dpar f^n}{\dpar\tb^j}(g(b))(w)+
$$
$$
+\frac{\dpar^2f^k}{\dpar\tb^i\dpar\tb^l}(g(b))\frac{\dpar g^l}{\dpar b^r}
\psi^k\phi^r(z)
\frac{\dpar^2f^n}
{\dpar\tb^j\dpar\tb^l}(g(b))\frac{\dpar g^l}{\dpar b^r}\psi^n\phi^r(w)
$$
When we compute each term using the Wick theorem, there 
appear single and double 
pairings. The part corresponding to the single pairings coincides 
with the expression of the bracket $[\dpar_{\tb^i},\dpar_{\tb^j}]$ 
in old coordinates $b^p$, so it vanishes. The "anomalous" part comes 
from the double pairings. 

One double pairing appears in the first term and is equal to  
$$
-\frac{\dpar}{\dpar b^n}\biggl[\frac{\dpar f^k}{\dpar\tb^i}(g(b))\biggr](z)
\frac{\dpar}{\dpar b^k}\biggl[\frac{\dpar f^n}{\dpar\tb^j}(g(b))\biggr](w)
\cdot\frac{1}{(z-w)^2}
$$
another one appears in the fourth term and equals 
$$
\delta_{km}\delta_{rn}
\frac{\dpar^2f^k}{\dpar\tb^i\dpar\tb^l}(g(b))\frac{\dpar g^l}{\dpar b^r}(z)
\frac{\dpar^2f^n}{\dpar\tb^j\dpar\tb^p}(g(b))\frac{\dpar g^p}{\dpar b^m}(w)
\cdot\frac{1}{(z-w)^2}=
$$
$$
=
\frac{\dpar^2f^k}{\dpar\tb^i\dpar\tb^l}(g(b))\frac{\dpar g^l}{\dpar b^n}(z)
\frac{\dpar^2f^n}{\dpar\tb^j\dpar\tb^p}(g(b))\frac{\dpar g^p}{\dpar b^k}(w)
\cdot\frac{1}{(z-w)^2}.
$$
We see that these terms cancel out. 

The remaining relations, (3.15c - f), contain only single pairings, 
and are easily checked. $\btu$  
  
Thus, for each automorphism $g=(g^1,\ldots,g^N)$ of $\BC[[b^1,\ldots,b^N]]$, 
(3.16a), the formulas (3.17) determine a morphism of vertex algebras 
$$
\tg:\ \hOmega_N\lra\hOmega_N
\eqno{(3.19)}
$$
More precisely, each element $c$ of $\hV_N$ is an (infinite) sum of  
finite products of $a^i_{-n}, b^j_{-m}, \psi^i_{-n}, \psi^j_{-m}$. 
We have 
$c=\bigl(c(z)\b1\bigr)(0)$. By definition, 
$$
\tg(c)=\biggl(\tg(c(z))\b1\biggr)(0)
\eqno{(3.20)}
$$
Thus, we have to define the field $\tg(c(z))$ for each $c\in\hV_N$. 
If $c$ is one of the generators $a^i_{-1}, b^i_0, \phi^i_0, \psi^i_{-1}$, 
we define $\tg(c(z))$ by formulas (3.17). We 
set 
$$
\tg(a^i_{-1-n}(z))=\dpar_z^{(n)}\tg(a^i_{-1}(z)),
\eqno{(3.21)}
$$
and the same with 
$b, \phi, \psi$, cf. (1.12). Finally, if $c=c^1c^2\cdot\ldots\cdot c^p$ where 
each $c^i$ is one of the letters $a, b, \phi$ or $\psi$, we set 
$$
\tg(c(z))=:\tg(c^1(z))\tg(c^2(z))\cdot\ldots\cdot\tg(c^p(z)):
\eqno{(3.22)}
$$
where the normal ordered product of two factors is defined by (1.2), 
and if $p>2$, we use the inductive formula (1.14). 

Equivalently, if $c^j=x^j_{k_j}$ where $x^j=a^i, b^i, \phi^i$ or $\psi^i$, 
we have
$$
\tg(c^1\cdot\ldots\cdot c^p\b1)=
\bigl[\tg(x^1(z))\bigr]_{k_1}\cdot\ldots\cdot
\bigl[\tg(x^p(z))\bigr]_{k_p}\b1
\eqno{(3.23)}
$$
(Here we use the following notation. If $a(z)=\sum_i a_iz^{-i-n}$ is 
a field corresponding to an element of conformal weight $n$, we denoted 
the Fourier mode $a_i$ by $a(z)_i$.)     

Let $G_N$ denote the group of automorphisms (3.16a). 

{\bf 3.8. Theorem.} {\it The assignment $g\mapsto\tg$ defines the 
group homomorphism $G_N\lra\Aut(\hOmega_N)$.} 

Let us consider two coordinate transformations, 
$'b^i=g_1^i(b)$, and $''b^i=g_2^i('b)$. Let $f_j$ denote the transformation 
inverse to $g_j$.  
We have to check that 
$$
\widetilde{g_2g_1}=\tg_2\tg_1
\eqno{(3.24)}
$$
By Theorem 4.5 from [K], 
it suffices to check this equality on the generators. Let us begin with  
$a^i$. The element $''a^i_{-1}\b1$ is expressed in the 
coordinates $'a$, etc., as follows  
$$
''a^i_{-1}\b1=\biggl(\ 'a^j_{-1}\frac{\dpar f_2^j}{\dpar''b^i}(g_2('b_0))-
\frac{\dpar^2f_2^k}{\dpar''b^i\dpar''b^l}(g_2('b_0))
\frac{\dpar g_2^l}{\dpar'b^r}('b_0)\ '\psi_{-1}^k\ '\phi_0^r\biggr)\b1
$$
Expressing it in the coordinates $a$, etc., we get the element
$$
\tg_2\tg_1(a^i)=
$$
$$
=\biggl[ a^p\frac{\dpar f^p_1}{\dpar'b^j}(g_1(b))(z)+
\frac{\dpar^2f_1^p}{\dpar'b^j\dpar'b^q}(g_1(b))\frac{\dpar g_1^q}
{\dpar b^s}\psi^p\phi^s(z)\biggr]_{-1}
\frac{\dpar f_2^j}{\dpar''b^i}(g_2g_1(b))(z)_0\b1-
$$
$$
-\frac{\dpar^2f_2^k}{\dpar''b^i\dpar''b^l}(g_2g_1(b))
\frac{\dpar g_2^l}{\dpar'b^r}(g_1(b))(z)_0
\biggl[\frac{\dpar f_1^p}{\dpar'b^k}(g_1(b))\psi^p(z)\biggr]_{-1}
\biggl[\frac{\dpar g_1^r}{\dpar b^q}\phi^q(z)\biggl]_0\b1
\eqno{(3.25)}
$$
(we have used (3.23)).  
Now, the action of our group $G_N$ on the classical de Rham complex 
is associative. It follows that the expression (3.25) is equal 
to $\widetilde{g_2g_1}(a^i)$ plus two anomalous terms: 
$$
a_0^p\biggl[\frac{\dpar f_1^p}{\dpar'b^j}(g_1(b))(z)\biggr]_{-1}
\frac{\dpar f_2^j}{\dpar''b^i}(g_2g_1(b_0))\b1=
\biggl[\frac{\dpar f_1^p}{\dpar'b^j}(g_1(b))(z)\biggr]_{-1}
\frac{\dpar}{\dpar b_0^p}\frac{\dpar f_2^j}{\dpar''b^i}(g_2g_1(b_0))\b1
$$
coming from the first summand, and
$$
-\frac{\dpar^2f_2^k}{\dpar''b^i\dpar''b^l}(g_2g_1(b_0))
\frac{\dpar g_2^l}{\dpar'b^r}(g_1(b_0))
\biggl[\frac{\dpar f_1^p}{\dpar'b^k}(g_1(b))(z)\biggr]_{-1}
\frac{\dpar g_1^r}{\dpar b^q}(b_0)\psi^p_0\phi^q_0\b1=
$$
$$
= - \frac{\dpar^2f_2^k}{\dpar''b^i\dpar''b^l}(g_2g_1(b_0))
\frac{\dpar g_2^l}{\dpar'b^r}(g_1(b_0))
\biggl[\frac{\dpar f_1^p}{\dpar'b^k}(g_1(b))(z)\biggr]_{-1}
\frac{\dpar g_1^r}{\dpar b^p}(b_0)\b1
$$
coming from the second one. 
One sees that these two terms cancel out, which proves (3.24) for $a^i$. 
For the generators $b, \psi$, and $\phi$, the anomaly does not appear at all. 
$\btu$  
 
{\bf 3.9.} Theorems 3.7 and 3.8 allow one to define the sheaf of 
conformal vertex algebras $\Omega^{ch}_X$ for each smooth manifold  
$X$ in an invariant way, by gluing the sheaves defined in 3.5. 
This can be done in each of 
the three settings: in algebraic, complex analytic or smooth one. 

In the complex analytic situation, we have our sheaves of vertex algebras 
over the coordinate charts, and the formulas (3.17) allow to glue these 
sheaves in a sheaf over $X$. 

In the algebraic situation, Theorem 3.8 ensures the existence 
of our sheaves by the standard arguments of "formal geometry" 
of Gelfand and Kazhdan, cf. [GK]. Consider the formal situation. 
By 3.8, the vertex algebra $\hOmega_N$ is a $G_N$-module. Therefore, 
the Lie algebra $W'_N=Lie(\hOmega_N)$ of formal vector fields vanishing 
at $0$ acts on $\hOmega_N$ by derivations. In fact, since the proof of Theorem 
3.8 
never uses the fact that the automorphisms in question preserve the origin, 
the infinitesimal version of formulas (3.17) shows that the entire 
algebra $W_N$ of formal vector fields operates on $\hOmega_N$.  
(Alternatively, this can be shown by the computation 
similar to the one from 5.1 (in the case of 
$\hOmega_N$ the anomaly vanishes!)). In other 
words, $\hOmega_N$ is a $(W_N,G_N)$-module (cf. [BS], [BFM]). 
Now, the standard result, [GK], says that such a module 
defines naturally a sheaf $\Omega^{ch}_X$ on each smooth algebraic variety $X$. 
They are sheaves of vertex algebras since $G_N$ (resp., $W_N$) 
acts by vertex algebra automorphisms (resp., derivations).  

A more direct construction of these sheaves is outlined in \S 6, see 6.10.

{\bf 3.10.} Consider the formal situation.  
Let us show that the algebra $\hOmega_N$ admits a canonical filtration whose 
graded factors 
are  standard tensor fields.
 
Placing the formulae (3.16) and (3.17), (3.23) on the desk 
next to each other, one realizes that the "symbols of fields" transform 
in the same way as the corresponding geometric quantities: functions, 
$1$-forms, and vector fields. To be more precise, introduce a filtration 
on $\hOmega_N$ as follows. 

The space $\hOmega_N$ is a free $\hA_N$-module with a base consisting of 
monomials in letters $a^i_n, \psi^i_n\ (n<0);\ b^i_m, \phi^i_m\ (m\leq 0)$. 
Define a partial ordering on this base by 

(a) $a>\phi,\ a>\psi,\ a>b,\ \psi>\phi,\ \psi>b,\ \phi>b$; $x^i_n>x^j_m$ 
if $n<m$, $x$ being $a, b, \phi$ or $\psi$; 

(b) extending this order to the whole set of monomials lexicographically. 

This partial order on the base naturally determines an increasing 
exhausting filtration 
on the spaces of fixed conformal weight,  
$$
F_0\hOmega_N^{(i)}\subset F_1\hOmega_N^{(i)}\subset\ldots
\eqno{(3.27)}
$$
For example, $F_0\hOmega^{(0)}_N=\hA_N,\ F_0\hOmega^{(-1)}_N=\oplus_{i=1}^N\ 
\hA_Nb_{-1}^i$, etc. A glance at (3.17), (3.23) shows that the corresponding 
graded object $\Gr^F_\bullet\Omega^{(i)}_N$ is a direct sum 
of symmetric powers 
of the tangent bundle, symmetric powers of the bundle of $1$-forms, 
and tensor products thereof. 
For example, the image of $b^i_{0}$ in $\Gr^F_\bullet\hOmega^{(0)}_N$ is 
a function, that of $b^i_{-n}\ (n>0)$ is a $1$-form, that of $a^i_{-n}$ is 
a vector field, etc. 

This filtration is stable under coordinate changes. Therefore 
all the sheaves $\Omega^{ch}_X$ acquire the natural filtration 
with graded factors being the bundles of tensor fields.

\bigskip

\bigskip
\centerline{\bf \S 4. Conformal and topological structure}
\bigskip

{\bf 4.1.} Let us return to the formal setting 3.6. Recall 
that we have in our vertex algebra $\hOmega_N$ the fields 
$L(z), J(z), Q(z)$ and $G(z)$, defined by the formulas (2.4), 
which make it a topological vertex algebra. 
Let us study the effect of the coordinate changes (3.16a) on these fields. 

Let us denote by $\tL(z)$, etc., the field $L(z)$, etc., written down 
using   
formulas (2.4), in terms of the tilded fields $\ta^i(z)$, etc., and then 
expressed in terms of the old fields $a^i(z)$, etc. 

{\bf 4.2. Theorem.} {\it We have
$$
\tL(z)=L(z)
\eqno{(4.1a)}
$$ 
$$
\tJ(z)=J(z)+\bigl(\Tr\log(\dpar g^i/\dpar b^j(z))\bigr)'
\eqno{(4.1b)}
$$
$$
\tQ(z)=Q(z)+\biggl(\frac{\dpar}{\dpar\tb^r}\bigl[
\Tr\log(\dpar f^i/\dpar\tb^j)\bigr]\tphi^r(z)\biggr)'
\eqno{(4.1c)}
$$
$$
\tG(z)=G(z)
\eqno{(4.1d)}
$$}

We have $J=\phi_0^i\psi_{-1}^i\b1$. Therefore, (cf. (3.23)), 
$$
\tJ=\biggl[\frac{\dpar g^i}{\dpar b^j}\phi^j(z)\biggr]_0
\biggl[\frac{\dpar f^k}{\dpar\tb^i}\psi^k(z)\biggr]_{-1}\b1=
J+\delta_{jk}\biggl[\frac{\dpar g^i}{\dpar b^j}(z)\biggr]_{-1}
\biggl[\frac{\dpar f^k}{\dpar\tb^i}(z)\biggr]_0\b1=
$$
$$
=J+\biggl(\frac{\dpar g^i}{\dpar b^j}(z)\biggl)_0'
\biggl(\frac{\dpar f^j}{\dpar\tb^i}(z)\biggr)_0\b1, 
$$
which implies (4.1b). We have $G=\psi^i_{-1}b^i_{-1}\b1$. Therefore, 
$$
\tG=\biggl[\frac{\dpar f^j}{\dpar\tb^i}\psi^j(z)\biggr]_{-1}
\bigl[g^i(b)(z)\bigr]_{-1}\b1=
\biggl[\frac{\dpar f^j}{\dpar\tb^i}\psi^j(z)\biggr]_{-1}
\biggl[\frac{\dpar g^i}{\dpar b^k} b^k(z)'\biggr]_0\b1=
$$
$$
=\delta_{jk}\psi^j_{-1}b^k_{-1}\b1=G
$$
This proves (4.1d). 
We have $Q=a^i_{-1}\phi_0^i\b1$. Therefore, 
$$
\tQ=\biggl[a^j\frac{\dpar f^j}{\dpar\tb^i}(z)-\frac{\dpar^2f^k}
{\dpar\tb^i\dpar\tb^l}\frac{\dpar g^l}{\dpar b^r}\psi^k\phi^r(z)\biggr]_{-1}
\biggl[\frac{\dpar g^i}{\dpar b^q}\phi^q(z)\biggr]_0\b1
$$
The classical terms: 
$$
\biggl[a^j_{-1}\biggl(\frac{\dpar f^j}{\dpar\tb^i}\biggr)_0-
\biggl(\frac{\dpar ^2f^k}{\dpar\tb^i\dpar\tb^l}\frac{\dpar g^l}{\dpar b^r}
\biggr)_0\psi^k_{-1}\phi^r_0\biggr]\biggl(\frac{\dpar g^i}{\dpar b^q}\phi^q
\biggl)_0\b1=Q, 
$$
since the second summand is zero, due to the anticommutation of 
$\psi^k_{-1}$ and $\phi^q_0$. Quantum corrections (anomalous terms): 
$$
\biggl[\biggl(\frac{\dpar f^j}{\dpar\tb^i}\biggr)_{-1}\biggl(\frac{\dpar^2g^i}
{\dpar b^j\dpar b^q}\biggr)_0\phi_0^r+
\biggl(\frac{\dpar^2f^k}{\dpar\tb^i\tb^l}\frac{\dpar g^l}{\dpar b^r}\phi^r
\biggr)_{-1}\biggl(\frac{\dpar g^i}{\dpar b^k}\biggr)_0\biggr]\b1=
$$ 
$$
=\biggl(\biggl[\biggl(\frac{\dpar f^j}{\dpar\tb^i}\biggr)_{-1}
\biggl(\frac{\dpar^2g^i}{\dpar b^j\dpar b^r}\biggr)_0+
\biggl(\frac{\dpar^2f^k}{\dpar\tb^i\dpar\tb^l}\frac{\dpar g^l}{\dpar b^r}
\biggr)_{-1}
\biggl(\frac{\dpar g^i}{\dpar b^k}\biggr)_0\biggr]\phi_0^r+
$$
$$
+\biggl(\frac{\dpar^2f^k}{\dpar\tb^i\dpar\tb^l}\frac{\dpar g^l}
{\dpar b^r}\frac{\dpar g^i}{\dpar b^k}\biggr)_0\phi^r_{-1}\biggr)\b1=
\biggl(\frac{\dpar f^k}{\dpar\tb^i\dpar\tb^l}\frac{\dpar g^l}
{\dpar b^r}\frac{\dpar g^i}{\dpar b^k}\phi^r\biggr)_{-1}\b1,
\eqno{(4.2)}
$$
since
$$
\biggl(\frac{\dpar^2f^k}{\dpar\tb^i\dpar\tb^l}\frac{\dpar g^l}
{\dpar b^r}\biggr)_0\biggl(\frac{\dpar g^i}{\dpar b^k}\biggr)_{-1}\b1=
-\biggl(\frac{\dpar f^k}{\dpar\tb^l}\frac{\dpar^2g^l}{\dpar b^t\dpar b^r}
\frac{\dpar f^t}{\dpar\tb^i}\biggr)_0\biggl(\frac{\dpar g^i}{\dpar b^k}
\biggr)_{-1}\b1=
$$
$$
=\biggl(\frac{\dpar f^k}{\dpar\tb^l}\frac{\dpar^2g^l}{\dpar b^t\dpar b^r}
\biggr)_0\biggl(\frac{\dpar f^t}{\dpar\tb^i}\biggr)_{-1}
\biggl(\frac{\dpar g^i}{\dpar b^k}\biggr)_0\b1=
\biggl(\frac{\dpar^2g^i}{\dpar b^t\dpar b^r}\biggr)_0
\biggl(\frac{\dpar f^t}{\dpar\tb^i}\biggr)_{-1}\b1
$$
Returning to (4.2), we have 
$$
\frac{\dpar^2 f^k}{\dpar\tb^i\dpar\tb^l}\frac{\dpar g^l}
{\dpar b^r}\frac{\dpar g^i}{\dpar b^k}\phi^r=
\frac{\dpar^2f^k}{\dpar\tb^l\dpar\tb^i}\frac{\dpar g^i}{\dpar g^k}\tphi^l=
\frac{\dpar}{\dpar\tb^l}\biggl(\Tr\log(\dpar f^i/\dpar\tb^j)\biggr)\tphi^l,
$$
which proves (4.1c). It follows from (4.1c) that the operator $Q_0$ is 
invariant. Hence, (4.1a) follows from (4.1d) and (2.2). $\btu$   

{\bf 4.3.} It follows from (4.1a) that for an arbitrary smooth manifold 
$X$, the field $L(z)$ is a well-defined global section of the sheaf 
$\End(\Omega_X^{ch})[[z,z^{-1}]]$, i.e. $\Omega^{ch}_X$ is canonically 
a sheaf of conformal vertex algebras. 

It follows from (4.1b) and (4.1c) that the Fourier modes $F=J_0$ ("fermionic 
charge") and $d_{DR}^{ch}=Q_0$ ("BRST charge") are well-defined 
endomorphisms of the sheaf $\Omega_X^{ch}$. Thus, 
$(\Omega^{ch}_X,d_{DR}^{ch})$ becomes a complex of sheaves, graded 
by $F$. This is a localization of definition 2.3.
 
{\bf 4.4. Theorem.} {\it For any smooth manifold $X$, the obvious 
embedding of complexes of sheaves 
$$
i:\ (\Omega_X,d_{DR})\lra (\Omega_X^{ch},d_{DR}^{ch})
\eqno{(4.3)}
$$
is a quasiisomorphism. This is true in algebraic, analytic and $C^{\infty}$  
settings.}

Indeed, the problem is local along $X$, and we are done by Theorem 2.4. 

{\bf 4.5.} If $X$ is Calabi-Yau, i.e. $c_1(\CT_X)=0$, then the fields 
$J(z)$ and $Q(z)$ are globally well defined, by (4.1b) and (4.1c). 
Here $\CT_X$ denotes the tangent bundle. Therefore, in this case 
the sheaf $\Omega_X^{ch}$ is canonically a sheaf of topological vertex 
algebras.

\bigskip
\centerline{\bf \S 5. Chiral structure sheaf}
\bigskip

{\it A. OBSTRUCTION}
\bigskip

{\bf 5.1.} Consider the formal setting 3.1, 3.2, 3.6. We have the 
vertex algebra $\hV_N$ of "chiral functions" over the formal disk 
$D_N=\Spf(\hA_N)$ where $\hA_N=\BC[[b^1,\ldots,b^N]]$. 
Let $W_N$ denote the Lie algebra of vector fields $f^i(b)\dpar_{b^i}$ 
on $D_N$. Let $\Omega^1(\hA_N)$ denote the module of one-forms 
$f^i(b)db^i$. The spaces $\hA_N$ and $\Omega^1(\hA_N)$ are naturally 
$W_N$-modules. Recall that the action of $W_N$ on $\Omega^1(\hA_N)$ 
is given by
$$
f^i\dpar_{b^i}\cdot g^jdb^j=f^i\dpar_{b^i}g^jdb^j+g^jdf^j
\eqno{(5.1)}
$$
The de Rham differential $d:\ \hA_N\lra\Omega^1(\hA_N)$ is compatible 
with the $W_N$-action.   

Let us define a map
$$
\pi:\ W_N\lra\End(\hV_N)
\eqno{(5.2)}
$$
For a vector field $\tau=f^i(b)\dpar_{b^i}$, let 
$\tau(z)$ denote the field $f^i(b)a^i(z)$ (of conformal 
weight $1$) of our vertex algebra 
$\hV_N$. Let $\pi(\tau)\in\End(\hV_N)$ denote the Fourier mode
$$
\pi(\tau):=\int\tau(z)=\tau(z)_0
\eqno{(5.3)}
$$
Note that by 1.7, the maps $\pi(\tau)$ are derivations of $\hV_N$.  

The mapping $\pi$ does {\it not} respect the Lie bracket. 
Let us compute the discrepancy. 
Let $\tau_1=f^i(b)\dpar_{b^i}, \tau_2=g^i(b)\dpar_{b^i}$ be two 
vector fields. We have the operator product
$$
\tau_1(z)\tau_2(w)=-\frac{\dpar_{b^j}f^i(b(z))\dpar_{b^i}g^j(b(w))}
{(z-w)^2}+
$$
$$
+
\frac{f^i(b(w))\dpar_{b^i}g^j(b(w))a^j(w)-g^j(b(w))\dpar_{b^j}f^i(b(w))a^i(w)}
{z-w}=
$$
$$
=-\frac{\dpar_{b^j}f^i(b(w))\dpar_{b^i}g^j(b(w))}
{(z-w)^2}+\frac{[\tau_1,\tau_2](w)}{z-w}-
\frac{\bigl[\dpar_{b^j}f^i(b(w))\bigr]'\dpar_{b^i}g^j(b(w))}{z-w}
\eqno{(5.4)}
$$
It follows that
$$
[\pi(\tau_1),\tau_2(w)]=[\tau_1,\tau_2](w)-
\bigl[\dpar_{b^j}f^i(b(w))\bigr]'\dpar_{b^i}g^j(b(w))
\eqno{(5.5)}
$$
In particular, 
$$
[\pi(\tau_1),\pi(\tau_2)]=\pi([\tau_1,\tau_2])-
\int\bigl(\dpar_{b^j}f^i(b(w))\bigr)'\dpar_{b^i}g^j(b(w))
\eqno{(5.6)}
$$

{\bf 5.2.} For $\omega=f^i(b)db^i\in\Omega^1(\hA_N)$, let us denote by 
$\omega(z)$ the field $f^i(b)b^i(z)'$ of our vertex algebra. Denote by 
$\pi(\omega)$ the Fourier mode $\omega(z)_0=\int\omega(z)$.

For $f=f(b)\in\hA^N$, let $f(z)$ denote the corresponding field of 
$\hV_N$. Its conformal weight is $0$, and $f(z)_0=f$. We have 
$$
df(z)=f(z)'
\eqno{(5.7)}
$$ 
Given $\tau=f^i\dpar_{b^i}\in W_N,\ \omega=g^jdb^j\in\Omega^1(\hA_N)$, 
we have the operator product
$$
\tau(z)\omega(w)=\frac{f^i(b(z))g^i(b(w))}{(z-w)^2}+
\frac{f^i(b(w))\dpar_{b^i}g^jb^j(w)'}{z-w}=
$$
$$
=\frac{f^i(b(w))g^i(b(w))}{(z-w)^2}+
\frac{f^i\dpar_{b^i}g^j(w)b^j(w)'+f^i(b(w))'g^i(b(w))}{z-w}
\eqno{(5.8)}
$$
It follows that
$$
[\pi(\tau),\omega(z)]=(\tau\omega)(z)
\eqno{(5.9)}
$$
and
$$
[\pi(\tau),\pi(\omega)]=\pi(\tau\omega)
\eqno{(5.10)}
$$
Let 
$\tW_N$ denote the linear subspace of $\End(\hV_N)$ generated by the 
Fourier modes $\pi(\tau)\ (\tau\in W_N)$ and $\pi(\omega)\ 
(\omega\in\Omega^1(\hA_N))$. Let $I_N\in\tW_N$ be the linear subspace 
generated by the Fourier modes $\pi(\omega)$. 

It follows from (5.7) that if $\omega$ is exact then $\pi(\omega)=0$. 
Thus, $\pi$ induces an epimorphic map 
$$
\Omega^1(\hA_N)/d\hA_N\lra I_N
\eqno{(5.11)}
$$

{\bf 5.3. Lemma.} {\it The map} (5.11) {\it is an isomorphism.}

This can be proved by writing down the Fourier mode as an infinite 
sum of monomials in $b^i_n$ and comparing the coefficients of like terms. 
In fact, a more general statement, namely, that $\int Q(z)=0$ if and 
only if $Q=\text{const}\cdot\b1$ or $Q(z)=P(z)'$ for some $P$, seems 
to be valid for a broad class of vertex algebras, cf. a similar 
statement in [FF3]. 

>From our point of view, this phenomenon has topological nature. It is 
amusing to exhibit an example of a vertex algebra, for which the lemma 
above is false. Namely, take $b^{-1}b_{-1}\in A_1[b^{-1}]$, see 3.4; then 
$(b^{-1}b_{-1})(z)=b(z)^{-1}b(z)'$ and $\int b(z)^{-1}b(z)'=0$, but 
$b(z)^{-1}b(z)'$ is not a total derivative. $\btu$    

{\bf 5.4.} Obviously, $\omega_1(z)\omega_2(w)=0$ for all 
$\omega_1, \omega_2\in\Omega^1(\hA_N)$. It follows from (5.6) and (5.10) 
that $\tW_N$ is a Lie algebra, $I_N$ is its abelian ideal, 
and we have the canonical extension 
$$
0\lra I_N\lra\tW_N\lra W_N\lra 0
\eqno{(5.12)}
$$
The action of $W_N$ on $I_N$ arising from this extension coincides with the 
canonical action of $W_N$ on $\Omega^1(\hA_N)/d\hA_N$, by (5.10). Note that 
we have defined this extension together with its splitting (5.2). It is 
given by the two-cocycle $c\in Z^2(W_N;\Omega^1(\hA_N)/d\hA_N)$ of 
$W_N$ with values in 
$\Omega^1(\hA_N)/d\hA_N$, read from (5.5), 
$$
c(f^i\dpar_{b^i},g^j\dpar_{b^j})= -\dpar_{b^i}g^jd(\dpar_{b^j}f^i)
(\text{mod}\ d\hA_N)
\eqno{(5.13)}
$$

{\bf 5.5.} Let us consider the truncated and shifted de Rham complex
$$
\Omega^\bullet:\ 0\lra\hA_N\lra\Omega^1(\hA_N)\lra 0 
\eqno{(5.14)}
$$
where we place $\Omega^1(\hA_N)$ in degree zero. It is a complex 
of $W_N$-modules. We have an obvious map of complexes of $W_N$-modules
$$
\Omega^\bullet\lra\Omega^1(\hA_N)/d\hA_N
\eqno{(5.15)}
$$
where the target is regarded as a complex sitting in degree zero. 

Let us write down a two-cocycle $\tc\in Z^2(W_N;\Omega^\bullet)$ which 
is mapped to $c$, (5.13), under the map (5.14). Such a cocycle is by 
definition a pair $(c^2,c^3)$ where $c^2\in C^2(W_N;\Omega^1(\hA_N))$ 
is a two-cochain, and $c^3\in C^3(W_N;\hA_N)$ is a three-cochain, such that 
$$
d_{Lie}(c^2)=d_{DR}(c^3)
\eqno{(5.16a)}
$$
$$ 
d_{Lie}(c^3)=0
\eqno{(5.16b)}
$$
Let us define 
$$
c^2(f^i\dpar_i,g^j\dpar_j)=\dpar_ig^jd(\dpar_jf^i)-
\dpar_jf^id(\dpar_ig^j)
\eqno{(5.17)}
$$
and
$$           
c^3(f^i\dpar_i,g^j\dpar_j,h^k\dpar_k)=\dpar_jf^i\dpar_kg^j\dpar_ih^k-
\dpar_kf^i\dpar_ig^j\dpar_jh^k
\eqno{(5.18)}
$$
We write for brevity $\dpar_i$ instead of $\dpar_{b^i}$. One checks the 
compatibilities (5.16) directly. Thus, we have defined $\tc$. One sees that 
$\tc$ is mapped to $-2c$ under (5.15). 

For $N=1$ the space $\Omega^1(\hA_N)/d\hA_N$ is trivial. This allows one  
to define the sheaf $\CO_X^{ch}$ for curves, acting as in 3.9, and 
starting from $\hV_N$ instead of $\hOmega_N$.  

Assume that $N>1$. 
Using the computations of Gelfand-Fuchs, cf. [F], Theorems 2.2.7 and 2.2.4,   
one can show that the map in cohomology
$$
H^2(W_N;\Omega^\bullet)\lra H^2(W_N;\Omega^1(\hA_N)/d\hA_N)
\eqno{(5.19)}
$$
induced by (5.15) 
is an isomorphism. We have the canonical short exact sequence
$$
0\lra H^2(W_N;\Omega^1(\hA_N))\lra H^2(W_N;\Omega^\bullet)\lra
H^3(W_N;\hA_N)\lra 0
\eqno{(5.20)}
$$
the left- and right-most terms being one-dimensional. Under 
the second map of this 
sequence, our cocycle $\tc$ is mapped to its second component $c^3$ 
which is a canonical representative of a generator of the space 
$H^3(W_N;\hA_N)$, cf. [F], Theorem $2.2.7'$ and Chapter 2, \S 1, no. 4. 
In particular, our cocycle $\tc$ is non-trivial. It follows that the cocycle 
$c$ defining the extension (5.12) is also non-trivial.  

What kind of an object does the cocycle $\tc$ define? Recall that a  
{\it homotopy Lie algebra} $L^\bullet$ is a complex of vector spaces 
equipped with a collection of brackets 
$$
[\ ,\ldots,\ ]_i:\ \Lambda^iL^\bullet\lra L^\bullet[-i+2],\ i\geq 2,
\eqno{(5.21)}
$$
satisfying certain compatibility conditions, cf. for example [HS], \S 4. 
In particular $[\ , ]_2$ is a skew symmetric map, satisfying the 
Jacobi identity up to the homotopy (given by the third bracket). 

Let us define a complex $L^\bullet$ as follows. Set 
$L^0=W_N\oplus\Omega^1(\hA_N)$, $L^{-1}=\hA_N$, the other components being 
trivial. The differential $L^{-1}\lra L^0$ is the composition of the de Rham 
differential and the obvious embedding $\Omega^1(\hA_N)\hra L^0$. 

Let the second bracket $[\ ,\ ]_2$ be given by the usual 
bracket of vector fields, and the action of vector fields on 
$\hA_N$ and $\Omega^1(\hA_N)$. Define the third bracket with the only 
nontrivial component being the three-cocycle $c^3$, (5.18). We set the other 
brackets equal to zero. This way we get a structure of a homotopy Lie algebra 
on $L^\bullet$. 

We have a canonical extension of homotopy Lie algebras 
$$
0\lra\Omega^\bullet\lra L^\bullet\lra W_N\lra 0
\eqno{(5.22)}
$$
Here $\Omega^\bullet$ is an abelian ideal in $L^\bullet$ (all brackets 
are zero). This is a refinement of extension (5.12).

\bigskip
{\it B. PROJECTIVE LINE}
\bigskip

{\bf 5.6.} Let $X$ be the projective line $\BP^1$. Let us fix a coordinate 
$b$ on $\BP^1$, and consider the open covering $X=U_0\cup U_1$ where 
$U_0=\Spec(\BC[b]),\ U_1=\Spec(\BC[b^{-1}])$. 

Consider the sheaves $\CO_{U_0}^{ch}$ on $U_0$ with coordinate $b$, and 
$\CO_{U_1}^{ch}$ on $U_1$ with coordinate $b^{-1}$, which 
were defined in 3.4. Let us glue them on the 
intersection $U_{01}=U_0\cap U_1$ using the transition function
$$
\tb(z)=b(z)^{-1}
\eqno{(5.23a)}
$$
$$
\ta(z)=b^2a(z)+2b(z)'
\eqno{(5.23b)}
$$  
In this way, we get the sheaf on the $X$, to be denoted $\CO_X^{ch}$. 

{\bf 5.7. Theorem.} {\it The space of global sections $\Gamma(X;\CO_X^{ch})$ 
admits a natural structure of an irreducible vacuum $\hsl_2$-module 
on the critical level.}

We have to compute $\CO_{U_0}^{ch}\cap\CO_{U_1}^{ch}$ where both 
$\CO_{U_i}^{ch}$ are regarded as subspaces of $\CO_{U_{01}}^{ch}$. It is 
the essence of the Wakimoto construction, [W], that the fields 
$a(z)b(z)^2+2b(z)', a(z)$ (resp., $\ta(z)\tb(z)^2+2\tb(z)', \ta(z)$) generate 
an $\hsl_2$-action on $\CO_{U_0}^{ch}$ (resp., on $\CO_{U_1}^{ch}$), and 
under this action, $\CO_{U_i}^{ch},\ i=1,2,$ become the restricted 
Wakimoto module with zero highest weight. (Restricted here means 
that the level is critical, and the Sugawara operators act by zero.) 
It follows from (5.23) that the $\hsl_2$-action comes from 
$\Gamma(X;\CO_X^{ch})$ and, therefore, $\CO_{U_0}^{ch}\cap\CO_{U_1}^{ch}$ 
is also an $\hsl_2$-module. It follows from [FF1] or [M] that each 
$\CO_{U_i}^{ch}$ contains a unique proper submodule which is isomorphic to 
the irreducible vacuum representation. To complete the proof, it remains 
to show that $\CO_{U_0}^{ch}\neq\CO_{U_0}^{ch}\cap\CO_{U_1}^{ch}$, 
and this is obvious. $\btu$

{\bf 5.8.} In fact, the first cohomology space $H^1(X;\CO_X^{ch})$ is 
also isomorphic to the same irreducible $\hsl_2$-module. 

To prove this, let us compute the Euler character
$$
ch(X;\CO_X^{ch})=\sum_{N=0}^\infty\ \chi(X;\CO_X^{ch(N)})\cdot q^N
$$
in two different ways. First, by definition
$$
ch(X;\CO_X^{ch})=ch(\Gamma(X;\CO_X^{ch}))-ch (H^1(X;\CO_X^{ch}))
$$
By Theorem 5.7 and [M], 
$$
ch (\Gamma(X;\CO_X^{ch}))=(1-q)^{-1}\prod_{N=1}^\infty\ (1-q^N)^{-2}
$$
On the other hand, formulas (5.23) imply that $\CO_X^{ch}$ carries 
a filtration $F$ such that the image of $a_{-n}$ (resp., $b_{-n}$) ($n\geq 1$) 
in $\Gr^F$ is a vector field (resp., a $1$-form). It follows that 
each monomial $a_{-n_1}\cdot\ldots\cdot a_{-n_r}b_{-m_1}\cdot\ldots\cdot
b_{-m_s}$ contributes $2s-2r+1$ in $\chi(X;\CO_X^{ch(N)})$, where 
$N=\sum n_i+\sum m_j$. Therefore, 
$$
ch(X;\CO_X^{ch})=\prod_{N=1}^\infty\ (1-q^N)^{-2},
$$
hence
$$
ch(H^1(X;\CO_X^{ch}))=q\cdot ch(\Gamma(X;\CO_X^{ch}))
$$
In other words, $H^1(X;\CO_X^{ch})$ has the same (up to the shift 
by $q$) character as $\Gamma(X;\CO_X^{ch})$. Again by [M], these 
two spaces are isomorphic as $\hfg$-modules. 

\bigskip    
  
{\it C. FLAG MANIFOLDS}

\bigskip

{\bf 5.9.} Let $G$ be a simple algebraic group, $B\subset G$ a Borel subgroup, 
$N\subset B$ the maximal nilpotent subgroup. The manifold 
$N$ is isomorphic to the affine space, and is a $(\fg,B)$-scheme, where $B$  
acts by conjugation. Consider the Heisenberg vertex algebra $V$ 
associated with the affine space $N$. According to [FF2], $V$ admits a 
structure of a $\hfg$-module (Wakimoto module); in particular, $V$ is a 
$(\fg,B)$-module. Note that $x\in\fg$ acts on $V$ as $\int X(z)$ for 
some $X\in V$.    

Consequently, considered as an affine space, $V$ 
admits a structure of a $(\fg,B)$-scheme.  
Let $M$ be the algebra of functions on $V$. Proceeding as in 3.9, 
with $K=B$, $\hX=G$, and $X=G/B$, we get the sheaf of ind-schemes 
$$
U\mapsto\Spec(H_{\nabla}(\Delta(M)))
$$
on $X$. The sheaf of its $\BC$-points  
is called the {\it chiral structure sheaf} of $X$ and denoted by $\CO_X^{ch}$. 


{\bf 5.10.} If $G=SL(n)$ then the sheaf $\CO_X^{ch}$ admits a more 
explicit construction, using charts and gluing functions. In this case 
$X=GL(n)/(B\times\BC^*)$. The Weyl group is identified with the symmetric 
group $S_n$ and can be realized as the subgroup of $GL(n)$ consisting of 
permutation matrices. One checks that in terms of the Lie algebra 
$\gl(n)$, the simple permutation $r_i$ (interchanging $i$ and $i+1$) 
can be written as follows
$$
r_i=\exp(\pi\sqrt{-1}E_{ii})\exp(E_{i+1,i})\exp(-E_{i,i+1})
\exp(E_{i+1,i})
\eqno{(5.24)}
$$
where $E_{ij}\ (1\leq i,j\leq n)$ form the standard base of $\gl(n)$. 

The manifold $X$ is covered by $|S|=n!$ charts, the chart associated with 
an element $w\in S_n$ being $U_w=wNw_0B$, where $N\subset B$ is the 
unipotent subgroup consisting of all upper-triangular matrices and 
$w_0\in S_n$ is the element of maximal length. Let us identify 
$U_w$ with $N$ using the bijection $n\mapsto wnw_0B$. Under this 
identification, if $x\in U_{w_1}\cap U_{w_2}$, then the change 
from the coordinates determined by $U_{w_1}$ to the ones determined by  
$U_{w_2}$, is given by $x\mapsto w_2^{-1}w_1x$. 

Each $U_w$ may be  identified with the affine space $\BC^{n(n-1)/2}$. 
To define $\CO_X^{ch}$, we declare  
that $\CO_X^{ch}\bigl|_{U_w}=\CO_{U_w}^{ch}$ where the last sheaf is defined 
in 3.4. Now we have to glue these sheaves over the pairwise intersections 
in a consistent manner. 

Let $V$ denote the vertex algebra $V_{n(n-1)/2}=
\Gamma(U_w;\CO_{U_w}^{ch})$. First, we extend the $\hsl(n)$-module 
structure on $V$ to a $\hgl(n)$-module structure. For that, in addition to 
the formulae in [FF1], p. 279, define 
$$
E_{ii}(z)=-\sum_{j>i}\ b^{ij}a^{ij}(z)+\sum_{j<i}\ b^{ji}a^{ji}(z)
\eqno{(5.25)}
$$
It is easily checked that in this way we indeed get an action of 
$\hgl(n)$ on $V$. 

For any $A\in\End(V)$, introduce the formal exponent
$\exp(tA):\ V\lra V\otimes\BC[[t]]$, 
$$
\exp(tA)(v)=\sum_{i=0}^\infty\ \frac{A^i(v)}{i!}t^i
\eqno{(5.26)}
$$
Working over $\BC[[t]]$, we can easily compose such maps. Motivated by 
(5.24), set 
$$
r_i(t)=\exp(t\pi\sqrt{-1}\int E_{ii}(z))\cdot
$$
$$
\cdot\exp(t\int E_{i+1,i}(z))
\exp(-t\int E_{i,i+1}(z))\exp(t\int E_{i+1,i}(z))
\eqno{(5.27)}
$$
Thus, $r_i(t)$ is a map $V\lra V[[z,z^{-1}]][[t]]$. Note that 
$V[[t]]$ is naturally a vertex algebra. It contains the vertex subalgebra 
"generated by functions, rational in $t$", that is, 
$$
V(t)=R(N\times\BC)\otimes_A V\subset V[[t]]
\eqno{(5.28)}
$$
Here $A=\Gamma(N;\CO_N)$, and $R(N\times\BC)$ denotes the ring 
of rational functions on $N\times\BC$, regular at $N\times\{0\}$. 

{\bf 5.11. Lemma.} (a) {\it $r_i(t)V\subset V(t)$. Further, $r_i(t)V$ is 
generated by functions well-defined for any value of $t$.} 

(b) {\it $r_i(1)$ is well defined (by} (a){\it ), and determines a map 
$$
r_i(1):\ \Gamma(U_w;\CO^{ch}_{U_w})\lra\Gamma(U_w\cap U_{r_iw},\CO^{ch}_{U_w})
\eqno{(5.29)}
$$}

As $\int X(z)$, $X\in V$, is a derivation, the map $r_i(t)$ is 
an embedding of vertex algebras. Therefore, it is enough to 
compute $r_i(t)$ on generators. 

There are two types of generators 
in $V$: (i) those coming from the subspace $W\subset V$, canonically 
isomorphic to $\fg=sl(n)$, such that the fields $X(z)\ (X\in W)$ generate 
the action of $\hfg$;  
(ii) those coming from $\BC[b_n^{ij}]\b1\subset V$. Obviously, the 
endomorphisms 
$\int E_{ij}(z)\ (i\neq j)$, act on $W$ as $E_{ij}$ on $sl(n)$; in particluar, 
this action is nilpotent, and $r_i(t)$ is polynomial in $t$. As for 
$\BC[b^{ij}_n]\b1$, this subspace is identified with the symmetric 
algebra $S^\bullet(\Omega^1(U_w))$. For any 
$X\in\gl(n)\subset\Gamma(X;\CT_X)$, the element $\int X(z)$ acts on this 
space as the Lie derivative along $X$. Consequently, 
$\exp(\int X(z))$ maps $S^\bullet(\Omega^1(U_w))$ into 
$S^\bullet(\Omega^1(\exp(X)\cdot U_w))$. $\btu$ 

Repeated application of this lemma gives for any $v=r_{i_1}\cdot\ldots r_{i_k}
\in S_n$ the map 
$$
v(1)=r_{i_1}\cdot\ldots\cdot r_{i_k}:\ \Gamma(U_w;\CO_{U_w}^{ch})\lra 
\Gamma(U_w\cap U_{v\cdot w};\CO_{U_w}^{ch})
\eqno{(5.30)}
$$
Finally, to complete our construction of $\CO_X^{ch}$, glue the sheaves 
$\CO_{U_{w_1}}$ and $\CO_{U_{w_2}}$ using the maps
$$
\Gamma(U_{w_1};\CO_{U_{w_1}}^{ch})\lra
\Gamma(U_{w_1}\cap U_{w_2};\CO_{U_{w_2}}^{ch})\lla
\Gamma(U_{w_2};\CO_{U_{w_2}}^{ch})
\eqno{(5.31)}
$$
where the first (resp. second) arrow is $w_2w_1^{-1}(1)$ (resp., 
$w_1w_2^{-1}(1)$). Since the gluing maps are induced by the action 
of $S_n$, they are transitive, and the sheaf $\CO_X^{ch}$ is 
well defined.

{\bf 5.12. Example.} Let $\fg=sl(2)$. We have $N=\BA^1\subset X=\BP^1$; 
$b$ is the coordinate on $N$ such that generators of $\gl(2)$ act as 
the following vector fields: 
$E_{21}\mapsto -\dpar_b$, $E_{12}\mapsto b^2\dpar_b$, 
$E_{11}\mapsto b\dpar_b$. 

One easily calculates that  
$$\exp(E_{21}):\ b\mapsto b-1;\ 
\exp(-E_{12}):\ b\mapsto b/(b+1);\ 
\exp(\pi\sqrt{-1}E_{11}):\ b\mapsto -b
$$ 
The simple reflection is 
$$
r=\exp(\pi\sqrt{-1}E_{11})\exp(E_{21})\exp(-E_{12})\exp(E_{21}):\ 
b\mapsto b^{-1}
\eqno{(5.32)}
$$
To do the chiral analogue of this computation, recall (cf. [W]) that 
$E_{21}(z)=-a(z),\ E_{12}(z)=b^2a(z)+2b(z)',\ E_{11}(z)=ba(z)$. 
A direct computation using Wick theorem yields 
$$
r(1)a_{-1}\b1=(b^2a_{-1}+2b_{-1})\b1
\eqno{(5.33)}
$$
which coincides with (5.23b). Computation of $r(1)b\b1$ does not differ 
from the classical one, see (5.32).            

{\bf 5.13. Theorem.} {\it The space $\Gamma(X;\CO_X^{ch})$ admits a 
natural structure of a $\hfg$-module at the critical 
level, such that it is a submodule of the restricted Wakimoto module, 
and its conformal weight zero component is $1$-dimensional.} 

Proof is not much different from that of Theorem 5.7. The space of sections 
$W=\Gamma(N;\CO_X^{ch})$ over the big cell is the restricted Wakimoto module 
with the zero highest weight. Restricted here again means that the 
center of $U(\hfg)_{loc}$ acts trivially.  
As the action of $B$ on $V$ preserves 
$\hfg\subset\End(V)$, $\hfg$ is spanned by the Fourier modes of fields 
associated with certain elements of $\Gamma(X;\CO_X^{ch})$. Therefore, 
$\Gamma(X;\CO_X^{ch})$ is a $\hfg$-module at the critical level, and there 
arises a $\hfg$-morphism
$$
\Gamma(X;\CO_X^{ch})\lra W
\eqno{(5.34)}
$$
By construction, the sheaf $\CO_X$ admits a canonical filtration 
whose associated quotients are free $\CO_X$-modules of finite rank. Therefore, 
the map (5.34) is an injection. By construction, the conformal 
weight zero component of $\Gamma(X;\CO_X^{ch})$ is equal to 
$\Gamma(X;\CO_X)=\BC$. $\btu$

This theorem is less precise than 5.7, the reason being that the 
representation-theoretic result of [FF1], [M] mentioned in 5.7 is only 
available in the $sl(2)$-case. However, Theorem 5.1 of [FF1] makes it 
plausible that this representation-theoretic statement carries over 
to any $\fg$, and hence that $\Gamma(X;\CO_X^{ch})$ is in fact the 
irreducible vacuum $\hfg$-module. 

{\bf 5.14.} Note that the whole sheaf $\CO_X^{ch}$ admits a structure of a 
$\hfg$-module at the critical level. 

{\bf 5.15. Localization for non-zero highest weight.} For any 
integral weight of $\fg$, there exists a twisted analogue of $\CO_X^{ch}$, 
to be denoted by $\CL_{\lambda}^{ch}$. Its construction repeats word 
for word 5.9, except that the action of $(\fg,B)$ on $V$ is to be 
twisted by $\lambda$, see [W, FF1]. If $\lambda$ is a regular dominant 
weight, then Theorems 5.7 and 5.13, along with their proofs, generalize 
in the obvious way. For example, the word "vacuum" in the formulation 
of Theorem 5.7 is to be replaced with "highest weight $\lambda$". 
Likewise, the claim that "conformal weight zero component is 
$1$-dimensional" should be replaced with the following: "conformal weight 
zero component is isomorphic to $\Gamma(X;\CL_{\lambda})$", where 
$\CL_{\lambda}$ is the standard line bundle on $X$. 

One may want to regard $\CL_{\lambda}^{ch}$ as a sheaf of modules 
over a sheaf of vertex algebras, just as $\CO_X^{ch}$ is a sheaf of modules 
over itself. From this point of view, $\CO_X^{ch}$ is a chiral 
analogue of the sheaf $\CD_X$ of differential operators on $X$. Thus, 
one cannot expect $\CL_{\lambda}^{ch}$ to be a module over 
$\CO_X^{ch}$. We believe (and have checked this for $\fg=sl(2)$) that 
one can define a sheaf of vertex algebras $\CO_{\lambda}^{ch}$ which 
acts on $\CL_{\lambda}^{ch}$ and is locally isomorphic to $\CO_X^{ch}$. 
Therefore, this sheaf can be regarded as a chiral partner of the 
sheaf $\CD_{\lambda}$ of twisted differential operators on $X$.

\bigskip
\centerline{\bf \S 6. Alternative construction}
\bigskip

In this section we will outline another construction of 
our vertex algebras, and some generalizations. The details 
will appear in a separate paper, see [MS].  

{\bf 6.1.} Recall (cf. [K], 4.9) that a {\it graded} vertex algebra is a 
pair $(V,H)$ where $V$ is a vertex algebra and $H:\ V\lra V$ is a 
diagonalizable linear operator ({\it Hamiltonian}) such that 
$$
[H,a(z)]=z\dpar_za(z)+(Ha)(z)
\eqno{(6.1.1)}
$$
for each $a\in V$. For example, a conformal vertex algebra is graded 
by the Hamiltonian $L_0$. The eigenvalues of $H$ are called 
conformal weights. By $V^{(\Delta)}$ we will denote the eigenspace 
of conformal weight $\Delta$. 
We have 
$$
a_{(n)}b\in V^{(\Delta+\Delta'-n-1)}\ \text{for\ } a\in V^{(\Delta)}, 
b\in V^{(\Delta')}
\eqno{(6.1.2)}
$$

Let us call a graded vertex algebra {\it restricted} 
if it has no negative integer conformal weights.
 
Let us fix a restricted vertex algebra $V$. To simplify the notations, 
we assume that $V$ is purely even. We would like to describe the structure 
which is induced on the space $V^{\leq 1}:=V^{(0)}\oplus V^{(1)}$ by 
the vertex algebra structure on 
$V$. We omit all computations; all the claims below are deduced directly 
from the Borcherds identity [K], Proposition 4.8 (b) and from 
{\it op. cit.} (4.2.3). 

{\bf 6.2.} The operation $a_{(-1)}b$ in $V$ will be denoted simply by $ab$. 

(a) The space $V^{(0)}$ is a commutative associative unital $\BC$-algebra 
with respect to the operation $ab$. This algebra will be denoted $A$. 
The unity is the vacuum, to be denoted by $\b1$.

The space $A$ acts by the left multiplication on $V^{(1)}$. However, 
this does not make $V^{(1)}$ an $A$-module: the multiplication by $A$ is 
not associative in general. 

We have the map $L_{-1}:\ A\lra V^{(1)}$. Let $\Omega\subset V^{(1)}$ 
denote the subspace spanned by the elements $a\dpar b,\ a,b\in A$. 
Thus, $L_{-1}$ induces the map
$$
d:\ A\lra\Omega
\eqno{(6.2.1)}
$$ 

(b) The left multiplication by $A$ makes $\Omega$ a left $A$-module. 
We have 
$$
a (db)=(db)a\ \ (a,b\in A)
\eqno{(6.2.2)}
$$

(c) The map $d$ is a derivation, i.e. 
$$
d(ab)=adb+bda
\eqno{(6.2.3)}
$$

Let us denote by $\CT$ the quotient space $V^{(1)}/\Omega$. 

(d) The left multiplication by $A$ makes $\CT$ a left $A$-module. 

Consider the operation
$$
_{(0)}:\ V^{(1)}\otimes V^{(1)}\lra V^{(1)}
\eqno{(6.2.4)}
$$

(e) The operation (6.6.4) induces a Lie bracket on $\CT$, to be denoted 
$[\ ,\ ]$.

Consider the operation 
$$
_{(0)}:\ V^{(1)}\otimes A\lra A
\eqno{(6.2.5)}
$$

(f) The operation (6.2.5) vanishes on the subspace $\Omega\otimes A$, 
and induces on $A$ a structure of a module over the Lie algebra $\CT$. 

This action will be denoted by $\tau(a)\ (a\in A,\ \tau\in \CT)$. 

(g) The Lie algebra $\CT$ acts on $A$ by derivations,  
$$
\tau(ab)=\tau(a)b+a\tau(b)\ 
\eqno{(6.2.6)}
$$

(h) We have 
$$
(a\tau)(b)=a\tau(b)
\eqno{(6.2.7)}
$$

The properties (d) --- (h) mean that $\CT$ is a {\it Lie algebroid} over $A$. 

(i) The operation (6.2.4) induces a structure of a module over the 
Lie algebra $\CT$ on the space $\Omega$.

This action will be denoted $\tau(\omega)$ or 
$\tau\omega\ (\tau\in \CT,\ \omega\in\Omega)$.  

(j) We have 
$$
\tau(a\omega)=a\tau(\omega)+\tau(a)\omega\ (a\in A,\ \tau\in\CT,\ 
\omega\in \Omega)
\eqno{(6.2.8)}
$$

(k) The differential $d:\ A\lra\Omega$ is compatible with the  
$\CT$-module structure.

It follows from (j) and (k) that

(l) we have
$$
\tau(a db)=\tau(a)db + a d(\tau(b))\ \ (\tau\in\CT,\ a,b\ \in A)
\eqno{(6.2.9)}
$$

Consider the operation
$$
_{(1)}:\ V^{(1)}\otimes V^{(1)}\lra A
\eqno{(6.2.10)}
$$

(m) The map (6.2.10) vanishes on the subspace 
$\Omega\otimes\Omega$. Therefore, it induces the pairing 
$$
\langle\ ,\ \rangle:\ \Omega\otimes\CT\oplus\CT\otimes\Omega\lra A
\eqno{(6.2.11)}
$$
This pairing is $A$-bilinear and symmetric. 
We have 
$$
\langle \tau, a db\rangle=a\tau(b)\ \ (\tau\in\CT,\ a,b\in A)
\eqno{(6.2.12)}
$$ 

(n) We have 
$$
(a\tau)(\omega)=a\tau(\omega)+\langle\tau,\omega\rangle da\ \ (a\in A,\ 
\tau\in\CT,\ \omega\in\Omega)
\eqno{(6.2.13)}
$$

{\bf 6.3.} Let us denote by $\hCT$ the space $V^{(1)}/dA$.  
The operation (6.6.4) induces a Lie bracket on the space $\hCT$. 
The subspace $\Omega/dA\subset\hCT$ is an abelian Lie ideal. 
The adjoint action of $\CT=\hCT/(\Omega/dA)$ coincides with 
the action defined by (i) and (l). 

Thus, we have an extension of Lie algebras
$$
0\lra \Omega/dA\lra \hCT\lra \CT\lra 0
\eqno{(6.3.1)}
$$
Note that this extension is not central in general.

{\bf 6.4.} 
Let us denote the space $V^{(1)}$ by $\tCT$. Thus, we have an exact sequence 
of vector spaces 
$$
0\lra\Omega\lra\tCT\lra\CT\lra 0
\eqno{(6.4.1)}
$$
Both arrows are compatible with the left multiplication by $A$. Let $\pi$ 
denote the projection $\pi:\ \tCT\lra\CT$. 

Let us define the "bracket" $[\ ,\ ]:\ \Lambda^2\tCT\lra\tCT$ by 
$$
[x,y]=\frac{1}{2}(x_{(0)}y-y_{(0)}x)\ \ (x,y\in\tCT)
\eqno{(6.4.2)}
$$
This bracket does not make $\tCT$ a Lie algebra: the Jacobi identity 
is in general violated. Set
$$
J(x,y,z)=[[x,y],z]+[[y,z],x]+[[z,x],y]\ \ (x,y,z\in \tCT)
\eqno{(6.4.3)}
$$
Consider the operation (6.2.10). 

(a) The operation (6.2.10) is symmetric. It will be denoted 
by $\langle x,y\rangle$. 

Let us define the map $I:\ \Lambda^3\tCT\lra A$ by 
$$
I(x,y,z)=\langle x,[y,z]\rangle+\langle y,[z,x]\rangle+
\langle z,[x,y]\rangle
\eqno{(6.4.4)}
$$

(b) We have 
$$
J(x,y,z)=\frac{1}{6}dI(x,y,z)
\eqno{(6.4.5)}
$$

{\bf 6.5.} Let us choose a splitting
$$
s:\ \CT\lra\tCT
\eqno{(6.5.1)}
$$
of the projection $\pi$. Let us define the map
$$
\langle\ ,\ \rangle=\langle\ ,\ \rangle_s:\ S^2\CT\lra A
\eqno{(6.5.2)}
$$
by
$$
\langle\tau_1,\tau_2\rangle=\langle s(\tau_1),s(\tau_2)\rangle
\eqno{(6.5.3)}
$$
(we put the lower index $_s$ in the notation if we want to stress 
the dependence on the splitting $s$).
Let us define the map
$$
c^2=c^2_s:\ \Lambda^2\CT\lra \Omega
\eqno{(6.5.4)}
$$
by
$$
c^2(\tau_1,\tau_2)=s([\tau_1,\tau_2])-[s(\tau_1),s(\tau_2)]
\eqno{(6.5.5)}
$$
Let us define the map $K:\ \Lambda^3\CT\lra A$ by
$$
K(\tau_1,\tau_2,\tau_3)=\langle s(\tau_1),s([\tau_2,\tau_3])\rangle+
\langle s(\tau_2),s([\tau_3,\tau_1])\rangle+
\langle s(\tau_3),s([\tau_1,\tau_2])\rangle\
\eqno{(6.5.6)}
$$
Let us define the map
$$
c^3=c^3_s:\ \Lambda^3\CT\lra A
$$
by
$$
c^3(\tau_1,\tau_2,\tau_3)=-\frac{1}{2}K(\tau_1,\tau_2,\tau_3) 
+\frac{1}{3}I(s(\tau_1),s(\tau_2),s(\tau_3))
\eqno{(6.5.7)}
$$
cf. (6.4.4). Let us regard $c^2$ (resp. $c^3$) as Lie algebra cochains 
living in $C^2(\CT;\Omega)$ (resp., in $C^3(\CT;A)$). 

(a) We have
$$
d_{Lie}(c^2)=dc^3
\eqno{(6.5.8)}
$$

(b) We have 
$$
d_{Lie}(c^3)=0
\eqno{(6.5.9)}
$$

The identities (a) and (b) mean that the pair $c=(c^2,c^3)$ form a 
$2$-cocycle of the Lie algebra $\CT$ with coefficients in the 
complex $A\lra\Omega$. 

(c) We have 
$$
\langle[\tau_1,\tau_2],\tau_3\rangle+\langle\tau_2,[\tau_1,\tau_3]\rangle=
\tau_1(\langle\tau_2,\tau_3\rangle)-
\frac{1}{2}\tau_2(\langle\tau_1,\tau_3\rangle)-
\frac{1}{2}\tau_3(\langle\tau_1,\tau_2\rangle)-
$$
$$
-\langle \tau_2,c^2(\tau_1,\tau_3)\rangle -
\langle\tau_3,c^2(\tau_1,\tau_2)\rangle
\eqno{(6.5.10)}
$$

Let us investigate the effect of the change of a splitting. 
Let $s':\ \CT\lra \tCT$ be another splitting of $\pi$. The difference 
$s'-s$ lands in $\Omega$; let us denote it 
$$
\omega=\omega_{s,s'}:\ \CT\lra\Omega
\eqno{(6.5.11)}
$$
We regard $\omega$ as a $1$-cochain of $\CT$ with coefficients in $\Omega$. 
Let us define a $2$-cochain $\alpha=\alpha_{s,s'}\in C^2(\CT;A)$ by
$$
\alpha(\tau_1,\tau_2)=\frac{1}{2}\bigl(\langle \omega(\tau_1),\tau_2\rangle-
\langle \tau_1,\omega(\tau_2)\rangle\bigr)
\eqno{(6.5.12)}
$$

(d) We have 
$$
c^2_s-c^2_{s'}=d_{Lie}(\omega)-d\alpha
\eqno{(6.5.13)}
$$

(e) We have 
$$
-d_{Lie}(\alpha)=c^3_s-c^3_{s'}
\eqno{(6.5.14)}
$$

The properties (d) and (e) mean that
$$
c_s-c_{s'}=d_{Lie}(\beta)
\eqno{(6.5.15)}
$$
where $\beta=\beta_{s,s'}:=(\omega,\alpha)\in C^1(\CT;A\lra\Omega)$. 

Therefore, we have assigned to our vertex algebra a canonically defined 
"characteristic class"
$$
c(V)=c(V^{\leq 1})\in H^2(\CT;A\lra\Omega)
\eqno{(6.5.16)}
$$
as the cohomology class of the cocycle $c_s$.

{\bf 6.6.} Let us introduce the mapping 
$$
\gamma=\gamma_s:\ A\otimes\CT\lra\Omega
\eqno{(6.6.1)}
$$
by
$$
\gamma(a,\tau)=s(a\tau)-a s(\tau)
\eqno{(6.6.2)}
$$

(a) We have 
$$
\gamma(ab,\tau)=\gamma(a,b\tau)+a\gamma(b,\tau)+
\tau(a)db+\tau(b)da
\eqno{(6.6.3)}
$$

(b) We have 
$$
\langle a\tau_1,\tau_2\rangle=a\langle\tau_1,\tau_2\rangle +
\langle \gamma(a,\tau_1),\tau_2\rangle -\tau_1\tau_2(a)
\eqno{(6.6.7)}
$$  

(c) We have
$$
c^2(a\tau_1,\tau_2)=ac^2(\tau_1,\tau_2)+\gamma(a,[\tau_1,\tau_2])-
\gamma(\tau_2(a),\tau_1)+\tau_2\gamma(a,\tau_1)-
$$
$$
-\frac{1}{2}\langle\tau_1,\tau_2\rangle da+
\frac{1}{2}d(\tau_1\tau_2(a))-
\frac{1}{2}d(\langle\tau_2,\gamma(a,\tau_1)\rangle)
\eqno{(6.6.8)}
$$
$(a\in A,\ \tau_i\in\CT)$.

(d) We have
$$
c^3(a\tau_1,\tau_2,\tau_3)=ac^3(\tau_1,\tau_2,\tau_3)+
\frac{1}{2}\tau_1[\tau_2,\tau_3](a)-
$$
$$
-\frac{1}{2}\bigl\{\langle\tau_2,\gamma(a,[\tau_3,\tau_1])\rangle-
\langle\tau_3,\gamma(a,[\tau_2,\tau_1])\rangle+
\langle\tau_2,\gamma(\tau_3(a),\tau_1)\rangle-
\langle\tau_3,\gamma(\tau_2(a),\tau_1)\rangle\bigr\}+
$$
$$
+\frac{1}{2}\bigl\{\langle\tau_2,\tau_3\gamma(a,\tau_1)\rangle-
\langle\tau_3,\tau_2\gamma(a,\tau_1)\rangle\bigr\}
-\frac{1}{2}\langle [\tau_2,\tau_3],\gamma(a,\tau_1)\rangle
\eqno{(6.6.9)}
$$

6.7. Let us call a {\bf prevertex algebra} the data (a) --- (f) below. 

(a) A commutative algebra $A$. 

(b) An $A$-module $\Omega$, together with an $A$-derivation 
$d:\ A\lra \Omega$. We assume that $\Omega$ is generated as a vector 
space by the elements $a db\ (a,b\in A)$, i.e. the canonical map 
$\Omega^1(A):=\Omega^1_{\BC}(A)\lra\Omega$ is surjective. 

(c) An $A$-Lie algebroid $\CT$. Define the action of $\CT$ on $\Omega$ by 
$$
\tau(a db)=\tau(a) db+ a d(\tau(b)), 
\eqno{(6.7.1)}
$$
cf. (6.2.9). We assume that this action is well defined. 
It follows that $d$ is compatible with the action of $\CT$.

We assume that the formula 
$$
\langle\tau,a db \rangle=a\tau(b)
\eqno{(6.7.2)}
$$
gives a well defined $A$-bilinear pairing $\CT\times\Omega\lra A$.    

(d) A $\BC$-bilinear mapping $\gamma:\ A\times\CT\lra\Omega$ satisfying 
6.6 (a).

(e) A $\BC$-bilinear symmetric mapping $\langle\ ,\ \rangle:\ 
\CT\times\CT\lra A$ satisfying 6.6 (b).

(f) A $\BC$-bilinear skew symmetric mapping $c^2:\ \CT\times\CT\lra \Omega$. 
This map should 
satisfy 6.5 (c), 6.6 (c), and the property (6.7.7) below. Let us define the map 
$$
[\ ,\ ]':= [\ ,\ ]-c^2:\ \Lambda^2\CT\lra \CT\oplus\Omega
\eqno{(6.7.3)}
$$
Let us define the map 
$$
c^3:=-\frac{1}{2}\tK+\frac{1}{3}\tI:\ \Lambda^3\CT\lra A
\eqno{(6.7.4)}
$$
where
$$
\tK(\tau_1,\tau_2,\tau_3)=\langle\tau_1,[\tau_2,\tau_3]\rangle+
\langle\tau_2,[\tau_3,\tau_1]\rangle+
\langle\tau_3,[\tau_1,\tau_2]\rangle
\eqno{(6.7.5)}
$$
and
$$
\tI(\tau_1,\tau_2,\tau_3)=\langle\tau_1,[\tau_2,\tau_3]'\rangle+
\langle\tau_2,[\tau_3,\tau_1]'\rangle+
\langle\tau_3,[\tau_1,\tau_2]'\rangle
\eqno{(6.7.6)}
$$
In the last formula, we imply the symmetric pairing $\langle\ ,\ \rangle$ to be 
extended to the whole space $\CT\oplus\Omega$ using (6.7.2), (e), and 
setting it equal to zero on $\Omega\times\Omega$.

Now, with $c^3$ defined above, we require that 
$$
d_{Lie}(c^2)=dc^3;\ \ d_{Lie}(c^3)=0
\eqno{(6.7.7)}
$$ 

Let us call a restricted vertex algebra $V$ {\it split} if it is  
given together with a splitting (6.5.1). We have constructed in 6.2 --- 6.6 
a functor 
$$
\CP:\ (Split\ restricted\ vertex\ algebras)\ \lra\  (Prevertex\ algebras) 
\eqno{(6.7.8)}
$$

{\bf 6.7.1. Claim.} {\it Functor $\CP$ admits a left adjoint, to be denoted 
$\CV$.}

In other words, given a prevertex algebra $P=(A,\Omega,\CT,\ldots)$, 
the corresponding vertex algebra $\CV(P)$ is defined by its 
universal property: to give a morphism of vertex algebras from $\CV(P)$ 
to an arbitrary split restricted vertex algebra $V'$ is the same as to 
give a morphism of prevertex algebras $P\lra \CP(V')$. 
(Morphisms of prevertex algebras are defined in the obvious manner.)

The construction of $V=\CV(P)$ goes in two steps. First, the components $V_0, 
V_1$ and 
the operations $_{(i)},\ i=-2,-1,0,1$ acting on them, are recovered by inverting 
the discussion 6.2 --- 6.6. For example, $V^{(0)}=A;\ V^{(1)}=\CT\oplus\Omega$; 
$$
\tau_{1(0)}\tau_2=[\tau_1,\tau_2]-c^2(\tau_1,\tau_2)+
d\langle\tau_1,\tau_2\rangle
\eqno{(6.7.9)}
$$
$$
\tau_{1(1)}\tau_2=\langle\tau_1,\tau_2\rangle
\eqno{(6.7.10)}
$$
$(\tau_i\in\CT)$, etc. Now, the components of weights $\geq 2$ are recovered 
by "bootstrap" from the universal property.  

Note that the set of conformal weights of $\CV(P)$ is equal to $\BZ_{\geq 0}$ if 
$\CV(P)\neq\BC$.  

{\bf 6.8. Example.} Let $\CT$ be a Lie algebra over $\BC$ equipped 
with an invariant bilinear form $\langle\ ,\ \rangle$.   
Set $A=\BC$, $\Omega=0,\ c^2=0, \gamma=0$. This defines a prevertex algebra $P$. 
Note that the component defined by the rule 6.7 (f) is not equal to zero, 
but is given by 
$$
c^3(\tau_1,\tau_2,\tau_3)=-\frac{1}{2}\langle\tau_1,[\tau_2,\tau_3]\rangle
\eqno{(6.8.1)}
$$ 
The vertex algebra $\CV(P)$ coincides with the vacuum (level $1$) representation 
of the 
affine Kac-Moody Lie algebra corresponding to $(\CT,\langle\ ,\ \rangle)$. 

{\bf 6.9. Example.} Let $A$ be a $\BC$-algebra; set 
$\Omega=\Omega^1_{\BC}(A)$. Let $\CT_0$ be an abelian Lie algebra over 
$\BC$ acting by derivations on $A$.  

Set $\CT=A\otimes_{\BC}\CT_0$. There is a unique Lie bracket on $\CT$ 
making it a Lie algebroid over $A$, 
$$
[a\tau_1,b\tau_2]=a\tau_1(b)\tau_2-b\tau_2(a)\tau_1\ \ 
(\tau_i\in\CT_0, a, b\in A)
\eqno{(6.9.1)}
$$ 

We set $\langle\tau_1,\tau_2\rangle=0;\ \gamma(a,\tau)=0; 
c^2(\tau_1,\tau_2)=0; c^3(\tau_1,\tau_2,\tau_3)=0$ for 
$a\in A,\ \tau_i\in\CT_0$. Then the formulas 6.6 (a) --- (d) define the 
unique extension of the operations $\gamma, \langle\ ,\ \rangle, 
c^2$ and $c^2$ to the whole space $\CT$.

Namely,
$$
\gamma(a,b\tau)=-\tau(a)db-\tau(b)da
\eqno{(6.9.2)}
$$
$$
\langle a\tau_1,b\tau_2\rangle=-a\tau_2\tau_1(b)-b\tau_1\tau_2(a)-
\tau_1(b)\tau_2(a)
\eqno{(6.9.3)}
$$
It is convenient to write down $c=(c^2,c^3)$ as a sum of a simpler 
cocycle and a coboundary,
$$
c^2(a\tau_1,b\tau_2)=\ 'c^2(a\tau_1,b\tau_2)+d\beta(a\tau_1,b\tau_2)
\eqno{(6.9.4)}
$$
where 
$$
'c^2(a\tau_1,b\tau_2)=\frac{1}{2}\bigl\{\tau_1(b)d\tau_2(a)-
\tau_2(a)d\tau_1(b)\bigr\}
\eqno{(6.9.5)}
$$
$$
'c^2(a\tau_1,b\tau_2)=\frac{1}{2}\bigl\{\tau_1(b)d\tau_2(a)-
\tau_2(a)d\tau_1(b)\bigr\}
\eqno{(6.9.5)}
$$
$$
\beta(a\tau_1,b\tau_2)=
\frac{1}{2}\bigl\{b\tau_1\tau_2(a)-a\tau_2\tau_1(b)\bigr\}
\eqno{(6.9.6)}
$$
and
$$
c^3(a\tau_1,b\tau_2,c\tau_3)=\ 'c^3(a\tau_1,b\tau_2,c\tau_3)+
d_{Lie}\beta(a\tau_1,b\tau_2,c\tau_3)
\eqno{(6.9.7)}
$$
where
$$
'c^3(a\tau_1,b\tau_2,c\tau_3)=\frac{1}{2}\bigl\{
\tau_1(b)\tau_2(c)\tau_3(a)-\tau_1(c)\tau_2(a)\tau_3(b)\}
\eqno{(6.9.8)}
$$ 
This gives the a prevertex algebra $P$ . 

Note that the cocycle $('c^2,\ 'c^3)$ coincides with (minus one half of) 
the cocycle (5.17-18) if $A$ is the polynomial ring.

For example, let $A$ be smooth, and assume that there exists a base 
$\{\tau_i\}$  
of the left $A$-module $\CT:=Der_{\BC}(A)$ consisting of commuting 
vector fields. Let $\CT_0\subset\CT$ be the $\BC$-vector space spanned by 
$\{\tau_i\}$. 
This gives a prevertex algebra $P$. The vertex algebra $A^{ch}:=\CV(P)$ may be 
called 
a "chiralization of $A$". This definition depends on the choice of $\{\tau_i\}$, 
and this is esssential; when we change the basis, we may get a non-isomorphic 
vertex algebra: here the "anomaly" appears. 

Specifying even more, let $A$ be a polynomial algebra $A_N$, cf. 3.1. Let 
$\tau_i=\dpar_{b^i}$. Then the vertex algebra $\CV(P)$ coincides with the 
Heisenberg vertex   
algebra $V_N$. 

If $A'$ is an arbitrary commutative $A$-algebra given together with 
an action of $\CT_0$ 
extending its action on $A$, then the base change  
$P_{A'}=(A',\CT_{A'}:=A'\otimes_A\CT,\ldots)$ has an obvious structure of 
a prevertex algebra. 
We have $\CV(P_{A'})=A'\otimes_A\CV(P)$. This explains the remark about 
the base change in 3.3. 

There exists a common generalization of the above two examples.  

{\bf 6.10.} All the above considerations have an obvious "super" 
($\BZ/(2)$-graded) version. Let us consider the super version of the 
Example 2.9. Let us start from the de Rham superalgebra $\Omega_A$ 
of differential forms over an arbitrary smooth algebra $A$. 

Let us assume that there exists an \'etale map $\Spec(A)\lra \BA^N$ 
given by coordinate functions $\{b^i\}$ (this is true maybe after 
some Zariski localization). Lifting the coordinate vector fields 
$a^i=\dpar_{b^i}$ to $A$, we get an abelian base  
in the Lie algebra $Der(A)$ which gives rise to an abelian base in Lie 
superalgebra $Der(\Omega_A)$. Now, we proceed as in 6.9 (in its super version),  
and get a vertex superalgebra $\Omega^{ch}=\CV(P)$. The calculation  
in the proof of Theorem 3.7 
shows that this vertex superalgebra does not depend on the choice 
of local \'etale coordinates. 

This is nothing but $\Gamma(X;\Omega^{ch}_X)$ for $X=\Spec(A)$. 
This may be viewed as an alternative 
(or a version of) construction of the chiral de Rham complex.    

{\bf 6.11. Chiral Weyl modules.} Let us return to the even situation 
again. Let $V$ be a restricted vertex algebra, 
let $\CP(V)=(A,\Omega,\CT,\ldots)$ be the corresponding 
prevertex algebra. Assume that the Lie 
algebra $\CT$ coincides with $Der(A)$. 

Let $\CM$ be a graded vertex module over $V$ (the definition of such an object  
is an obvious modification of the definition of a graded 
vertex algebra). Assume that $\CM$ is {\it restricted}, i.e. has no negative 
integer conformal weights. 
Consider the weight zero component $M=\CM^{(0)}$. 
Then the operations $am=a_{(-1)}m\ (a\in A, m\in M)$ and 
$\tau m=x_{(0)}m\ (\tau\in\CT, x\in \tCT$ is any representative 
of $\tau$, $m\in M$) makes $M$ a left $\CD$-module over $A$. 

This way we get a functor
$$
(Restricted\ V-modules)\lra (D_A-modules)
\eqno{(6.11.1)}
$$
This functor admits a left adjoint 
$$
\CW:\ (D_A-modules) \lra (Restricted\ V-modules)
\eqno{(6.11.2)}
$$

For a $\CD$-module $M$, the vertex module $\CW(M)$ is called 
the {\it chiral Weyl module} corresponding to $M$.
 
It is natural to hope this construction applied to flag spaces $G/B$ gives 
a functor from $\CD$-modules over $G/B$ to modules over 
$\CO^{ch}_{G/B}$ which corresponds to the Weyl module construction 
in the language of representations. 

A similar construction gives for an arbitrary smooth manifold $X$ 
a functor 
$$
\Omega^{ch}:\ (\CD_X-mod)\lra (\Omega^{ch}_X-mod)
\eqno{(6.11.3)}
$$
called the {\it chiral de Rham complex of a $\CD$-module}.

\bigskip
\centerline{\bf References}
\bigskip

[BD1] A.~Beilinson, V.~Drinfeld, Chiral algebras I, Preprint. 

[BD2] A.~Beilinson, V.~Drinfeld, Quantization of Hitchin's integrable 
system and Hecke eigensheaves, Preprint, 1997.  

[BFM] A.~Beilinson, B.~Feigin, B.~Mazur, Introduction to algebraic 
field theory on curves, Preprint. 

[BS] A.~Beilinson, V.~Schechtman, Determinant bundles and Virasoro 
algebras, {\it Comm. Math.Phys.} {\bf 118} (1988), 651-701.

[B] L.A.~Borisov, Vertex Algebras and Mirror Symmetry, math.AG/9809094.  



[FF1] B.~Feigin, E.~Frenkel, Representations of affine Kac-Moody algebras 
and bosonization, in: V.~Knizhnik Memorial Volume, L.~Brink, D.~Friedan, 
A.M.~Polyakov (Eds.), 271-316, World Scientific, Singapore, 1990. 

[FF2] B.~Feigin, E.~Frenkel, Affine Kac-Moody algebras and semi-infinite 
flag manifolds, {\it Comm. Math. Phys.} {\bf 128} (1990), 161-189. 

[FF3] B.~Feigin, E.~Frenkel, Affine Kac-Moody algebras at the critical 
level and Gelfand-Dikii algebras, Infinite Analysis, Part A, B 
(Kyoto, 1991), 197-215, {\it Adv. Ser. Math. Phys.} {\bf 16}, 
World Sci. Publishing, River Edge, NJ, 1992.  

[F] D.B.~Fuchs, Cohomology of infinite-dimensional Lie algebras, 
Contemp. Sov. Math., Consultants Bureau, New York, 1986.

[GK] I.M.~Gelfand, D.A.~Kazhdan, Some problems of differential 
geometry and the calculation of cohomologies of Lie algebras of vector 
fields, {\it Soviet Math. Dokl.}, {\bf 12} (1971), No. 5, 1367-1370. 

[HS] V.~Hinich, V.~Schechtman, Homotopy Lie algebras, I.M.~Gelfand 
Seminar, {\it Adv. in Sov. Math.}, {\bf 16} (1993), Part 2, 1-28.       

[K] V.~Kac, Vertex algebras for beginners, University Lecture Series, 
{\bf 10}, American Mathematical Society, Providence, RI, 1997. 

 
[LVW] W.~Lerche, C.~Vafa, P.~Warner, Chiral rings in $N=2$ superconformal 
theories, {\it Nucl. Phys.} {\bf B324} (1989), 427-474.

[LZ] B.~Lian, G.~Zuckerman, New perspectives on the 
BRST-algebraic structure of string theory, {\it Commun. 
Math. Phys.} {\bf 154} (1993), 613-646; hep-th/9211072. 

[M] F.~Malikov, Verma modules over Kac-Moody algebras of rank $2$, 
{\it Leningrad Math. J.} {\bf 2}, No. 2 (1990), 269-286. 

[MS] F.~Malikov, V.~Schechtman, Chiral de Rham complex. II, 
D.B.~Fuchs' 60-th Anniversary Volume, 1999, to appear; 
math.AG/9901065.  

[W] M.~Wakimoto, Fock representations of the affine Lie algebra $A_1^{(1)}$, 
{\it Comm. Math. Phys.} {\bf 104} (1986), 604-609.   



\bigskip

F.M.: Max-Planck-Institut f\"ur Mathematik, 
Gottfried-Claren-Stra\ss e 26, 53225 Bonn, 
Germany;\ Department of Mathematics, University of 
Southern California, Los Angeles, CA 90089, USA;\  malikov\@mpim-bonn.mpg.de,\  
fmalikov\@mathj.usc.edu 

V.S.: Department of Mathematics, University of Glasgow, 15 University 
Gardens, Glasgow G12 8QW, United Kingdom;\ vs\@maths.gla.ac.uk

A.V.: Department of Mathematics, New Mexico 
State University, Las Cruces, NM 88003-8001, USA;\  
vaintrob\@nmsu.edu

\enddocument